\documentclass[a4paper,10pt]{article}

\usepackage[english]{babel}
\usepackage{amsmath}
\usepackage{siunitx}
\usepackage[a4paper,left=3cm,right=3cm,top=2.5cm,bottom=2.5cm]{geometry}
\usepackage{bm}
\usepackage{authblk}
\usepackage{siunitx}
\usepackage{graphicx}
\usepackage{float}
\usepackage{comment}
\usepackage{subfig}
\usepackage{mathtools}
\usepackage{amssymb}
\usepackage{amsthm}
\usepackage{multirow}
\usepackage{cite}
\usepackage{breakcites}
\usepackage[colorinlistoftodos]{todonotes}
\usepackage{tabularx}
\usepackage{longtable}
\usepackage[ruled, vlined, linesnumbered]{algorithm2e}
\usepackage{hyperref}
\usepackage{tikz}
\hypersetup{
colorlinks=true,
linkcolor=blue,
urlcolor=red,
pdftitle={Valentin1stPaper},
}

\usepackage{ifpdf}

\ifpdf
\usepackage{pdfsync}
\fi

\usepackage{cleveref}
\usepackage{algpseudocode}
\DeclareGraphicsExtensions{.pdf,.png}
\usepackage{caption}
\usepackage{subcaption}
\newcommand{\RA}[1]{{#1} }
\newcommand{\RB}[1]{{#1}}
\graphicspath{ {./images/} }

\title{A reduced-order model for segregated fluid-structure interaction  solvers based on an ALE approach}
\author[1]{Valentin Nkana Ngan\footnote{vkanang@sissa.it}}
\author[2]{Giovanni Stabile\footnote{giovanni.stabile@santannapisa.it}}
\author[3]{Andrea Mola\footnote{andrea.mola@imtlucca.it}}
\author[1]{Gianluigi Rozza\footnote{grozza@sissa.it}}

\affil[1]{Mathematics Area, mathLab, SISSA, via Bonomea 265, I-34136 Trieste, Italy}
\affil[2]{The Biorobotics Institute, Sant'Anna School of Advacended Studies, V.le R. Piaggio 34, 56025, Pontedera, Pisa - Italy}
\affil[3]{MUSAM Continuum Mechanics Laboratory, Scuola IMT Alti Studi Lucca - Piazza S. Ponziano, 6-55100 Lucca, LU, Italy }
\date{\today} 

\begin{document}
\tableofcontents
\maketitle
\begin{abstract}
This article presents a Galerkin projection-based reduced-order modeling (ROM) approach for segregated fluid–structure interaction (FSI) problems, formulated within an Arbitrary Lagrangian–Eulerian (ALE) framework at low Reynolds numbers using the Finite Volume Method (FVM). The ROM is constructed using Proper Orthogonal Decomposition (POD) and incorporates a data-driven technique that combines classical Galerkin projection with radial basis function (RBF) networks.

The results demonstrate the numerical {\RA{stability and accuracy of the proposed method relative to the high-fidelity model.}}The ROM successfully captures transient flow fields and, importantly, the forces acting on the moving structure without exhibiting unphysical growth or divergence over time. This is further supported by the bounded evolution of error metrics and physical observables, which remain consistent with the full-order simulations throughout the prediction horizon.

The method's effectiveness is validated through a benchmark vortex-induced vibration (VIV) case involving a circular cylinder at Reynolds number $Re = 200$. The hybrid ROM approach yields an accurate and efficient tool for solving FSI problems involving mesh motion.
\end{abstract}

\section{Motivation and state-of-the-art}
Understanding fluid-structure interaction (FSI) is essential for addressing numerous important engineering problems. The interaction between fluids and movable or deformable structures is of both historical and practical significance, playing a critical role in the design of various systems
\cite{Falaize2019, dipod}. In systems engineering, the behavior of bluff bodies---such as cylinders, prisms, and square cross-sections---is particularly relevant in fields like aerodynamics, automative engineering, and civil infracstructure \cite{buresti2022bluff}.
These structures, which are commonly subjected to flow---wheter air or water---are found in bridges, chimneys, marine cables, and pipelines. Accurately predicting their dynamic response is essential because they are prone to significant vibrations under fluid loading. Since the middle 20th century, this area has attracted considerable attention from researchers \cite{kara2015calculation}.

To simulate these systems numerically, flow past a circular cylinder has become a canonical benchmark problem in fluid dynamics. This configuration captures many fundamental phenomena and serves as a proxy for understanding more complex real-world applications\cite{Placzek2009, Placzek2007}, often described as a "kaleidoscope of challenging fluid phenomena" \cite{deane1991low}.

However, simulating FSI problems remains computationally expensive, even on modern high-performance computing platforms, due to the cost of storing, handling, and processing large-scale nonlinear data. To address this, significant research has been devoted to the development of reduced-order models (ROMs) \cite{garelli2011fluid}.

A ROM is a simplified mathematical representation of a physical system, derived from computational or experimental data, that retains the essential dynamics of the original model while significantly reducing its degrees of freedom. These models offer a cost-effective and efficient means of exploring  complex physical problems. The primary motivation behind ROMs is to achieve accurate physical insights with reduced computational overhead  \cite{mifsud2008reduced}.
While various ROM strategies exist, most share the common goal of identifying dominant spatial and temporal features in the flow field, often based on high-fidelity data from full-order simulations governed by the Navier–Stokes equations.

Among model reduction techniques, the POD-Galerkin (POD-G) method is widely used due to its success in various applications \cite{anttonen2001techniques}. Many researchers have explored ROMs for FSI problems, leading to a rich body of literature.  For instance, \cite{Nonino2021}, presents a combination of the Reduced Basis Method (RBM) with monolithic and partitioned FSI approaches, tested on the Turek–Hron benchmark with {\RB {$Re=100$}}. Optimization-based domain decomposition methods for incompressible Navier–Stokes equations were studied in \cite{prusak2022optimisationbased}, demonstrating significant computational gains.  Further, 
\cite{xu2021global} extended POD-Galerkin to domains with moving solid boundaries, while 
 \cite{liberge2010reduced, Falaize2019}, used POD with immersed interface methods to simulate flow over oscillating cylinders.
Some works, such as \cite{shinde2019galerkin}, avoid Galerkin projection due to numerical stability issues, using separate ROMs for fluid and structure domains instead. Others, like \cite{xiao2016non}, propose non-intrusive ROMs combining POD with RBF interpolation for unstructured finite element meshes. In \cite{shinde2016modelling}, decoupled ROMs were proposed, using singular value decomposition (SVD) to construct structural basis functions. The work by Liberge and Hamdouni  \cite{liberge2010reducedcyl}, building on their earlier 1D Burgers equation study \cite{liberge2007proper},  developed global POD modes for FSI problems involving spring-mounted cylinders.

Several hybrid ROM–machine learning approaches have also emerged
\cite{gupta2022hybrid, patel2021reduced, whisenant2020galerkin, miyanawala2019hybrid, miyanawala2019decomposition} to address the stability challenges inherent in Galerkin projection. Notably, most of these studies rely on the Finite Element Method (FEM) for the full-order model. 
However, recent efforts have introduced ROMs within the Finite Volume Method framework
\cite{hijazi2020data, kelbij2019novel, stabile2020efficient, ivagnes2023hybrid, dao2021projection}.
For example, \cite{dao2021projection} proposed a ROM for transient, multi-object cross-flow simulations using moving domains and immersed boundary methods and \cite{hijazi2020data} studied a reduced-order model based on the 
merging/combining projection-based techniques with data-driven reduction strategies on the a fixed mesh.

This work distinguishes itself from previous studies by incorporating mesh motion within the Arbitrary Lagrangian–Eulerian (ALE) framework, as a strong coupling between the fluid and structural domains. It further introduces a hybrid approach, where the primary partial differential equations (PDEs) are addressed using a standard POD–Galerkin projection, while Radial Basis Function (RBF) networks are employed as a data-driven method to interpolate mesh node displacements in a point cloud setting.
This choice is based on both theoretical studies and practical considerations. 
The practical aspects are related to the idea of generating an approach that could be applied to construct a reduced-order model independent of the mesh motion technique used at the full-order level.
Secondly, despite a large amount of theoretical work behind the mesh motion model, there are still several empirical coefficients, making the overall formulation less rigorous in terms of physical principles. 
These considerations have been used to propose a reduced-order model that could be applied to any mesh motion model that exploits a projection-based technique for mass and momentum equations and a data-driven approach for predicting the grid node displacement field. This work considers a  cross-flow of Reynolds number $Re = 200$. Although the test case is relatively simplified compared to realistic problems,  this study is among the few attempts to extend ROM's capability to model moving objects in the Finite Volume platform and capture the forces acting on the object. The development of ROM for simulations of high Re flow is beyond the focus of this paper. It is currently conducted in a separate analysis by the authors on different benchmark problems.

The manuscript is structured as follows: Section \ref{genmathform} begins with the formulation of the structure dynamic in  Subsection \ref{structmotion}. 
In the next three subsections, the mathematical formulation of the fluid's motion in the ALE setting is presented in see Subsection \ref{fluidmotion}, followed by the coupling strategy at the interface in Subsection \ref{couplingcond} and the discussion of mesh motion strategies in Subsection  \ref{motionstrategies}. 
Section \ref{NumDisFOMROM} starts by addressing the numerical discretization of the full-order model (see Subsection \ref{standardfv}), then Subsection \ref{rompro} addresses the 
reduced-order model concept; follows with the discussion of Proper Orthogonal Decomposition (POD) for laminar incompressible flows in  Subsection \ref{sectredpimple} by insisting on its main properties in terms of model reduction while Subsection \ref{pod-i} formally introduces the POD-RBF with radial basis functions networks. Section \ref{resdiscuss} discusses and presents the numerical results. {\RB {Section \ref{conclusion} concludes the manuscript and outlines potential directions for future research.}}
\section{Mathematical formulation of the fluid-structure interaction problem}\label{genmathform}
This section presents the mathematical formulation of the  fluid-structure interaction model. In this work, the  following assumptions are considered :
\begin{itemize}
    \item The fluid is viscous, incompressible and Newtonian, 
     \item A 2D simulation is carried out,
    \item The cylindrical structure is rigid. Its elastic connection to the ground is represented by a spring and a damper,  
    \item The cylinder is free  to oscillate in the vertical direction.
\end{itemize}
Based on these assumptions, the  present {\RB {section}} briefly describes the full-order model that simulates the coupled dynamics of both the fluid and structure. The fluid motion is governed by the Navier-Stokes equations in an arbitrary Lagrangian Eulerian (ALE) framework, while the structural motion is modeled using rigid-body dynamics with appropriate coupling conditions at the fluid-structure interface.

\subsection{Structural Dynamics Formulation}
\label{structmotion}

{\RA{
Vortex shedding flow past elongated crossflow cylinders is for the most part a two dimensional phenomenon. In fact, tip effects are not affecting the flow in the vast majority of the cylinder sections, which are located a considerable number of diameters away from the tips. As such, 2D simulations on central sections of the cylinder represent a good approximation of the local 3D flow field. A similar approximation is particularly suitable for the simulation of the flow past electricity or bridge cables, and offshore platforms risers.}  } In the 2D set-up we then consider, the structure of the cylinder circular cross section is considered rigid. In fact, the cylinder round section is not experiencing significant in plane deformations associated to the fluid dynamic forces. Instead, these forces, which are approximately constant along the $z$ direction indicated in Fig. \ref{fig:cross-flow-viv} (a), induce low curvature bending deformations along the cylinder length, which might however result in considerable displacement of the cylinder within the 2D sectional plane. Given these considerations, the elastic forces associated with the three-dimensional bending deformation are recovered in the sectional 2D model by allowing for free cylinder translations, and connecting the cylinder to the ground with a spring and a damper. To further simplify the problem, in this work the cylinder is constrained to translate only along the vertical direction ($y$). In the literature, this is a rather common additional assumption, and it is also consistent with experimental campaigns carried out to characterize vortex induced vibration (see for instance \cite{Khalak1997FluidFA}).
The vertical motion $y_C$(t) of the cylinder is described by the following second-order differential equation,
\begin{align}\label{ode}
	\ddot{y}^C + 2\zeta \omega_n\dot{y}^C + \omega^2_n y^C = \frac{F_y}{m},
\end{align}
where,
\begin{itemize}
\item $\omega_n = \sqrt{k/m} = 2\pi f_n$ is the natural pulsation of
the system
\item $k$ is the spring's stiffness
\item $m$ is the mass of the rigid body
\item  $\zeta = \dfrac{c}{2m\omega_n}$ is the damping ratio, with $c$ being  the damping coefficient
\item $F_y$ is the lift force per unit length in the transverse direction
\item $y^C$ is the vertical displacement of the cylinder.
\end{itemize}
{\RA{
Referring to equation \eqref{ode}, it is important to remark that lift force $F_y$ is computed based on the fluid dynamic problem solution. It acts therefore as the coupling term from the fluid domain, driving the cylinder's motion. In addition, we must point out that the scalar equation governing the structural problem results, at the discrete level, only in a single degree of freedom. As will be seen in Section \ref{pod-i}
, such a computationally inexpensive structural problem will be used not only in the full order FSI problem. It will also be used in the reduced FSI problem, coupled with the reduced fluid model.}
}
\subsection{Fluid Dynamics Formulation}\label{fluidmotion}

The motion of the cylindrical structure in the transverse direction has the direct consequence that the fluid dynamic domain $\Omega=\widehat{\Omega}(t)$ a function of time. It in fact deforms to cope with the rigid displacements of one of its boundaries --- the one associated with the cylinder. We then make use of the Arbitrary Lagrangian-Eulerian formulation \cite{doneaALE2004} of the Navier--Stokes Equations to deal with the motion deformation of the computational domain.

The momentum equation in the ALE framework is written as follows:
\begin{align}\label{momentum}
\frac{\delta \boldsymbol{u} }{\delta t} + \nabla \cdot [\boldsymbol{u}\otimes(\boldsymbol{u} - \boldsymbol{u}^g)] -\nabla \cdot(\mu\nabla \boldsymbol{u}) = -\nabla p.
\end{align}
Where,
\begin{itemize}
    \item  $\boldsymbol{u}$ is the velocity
    \item $\boldsymbol{u}^g$ the grid velocity
    \item $p$ being the pressure
    \item $\mu$ the dynamic viscosity
    \item $\frac{\delta }{\delta t}$  is the ALE time derivative defined as,
    \begin{align}
    \frac{\delta }{\delta t} = \frac{\partial }{\partial t} + \boldsymbol{u}^g\nabla.
\end{align}
This derivative measures the rate of change of a quantity at grid-moving points.
\end{itemize} 
{\RA{
We point out that, in the ALE framework, the grid velocity field $\boldsymbol{u}^g$ is obtained as the time derivative of the grid displacement field $\boldsymbol{d}^g$, 

\begin{equation}
    \boldsymbol{u}^g = \boldsymbol{\dot{d}}^g.
\end{equation}
As will be discussed in Section~\ref{motionstrategies}, $\boldsymbol{d}^g$ is computed at each time step by means of a mesh morphing algorithm which smoothly extends the motion of the solid interface boundary to the internal nodes of the fluid grid. 
 }}
 
The incompressibility condition remains:
\begin{align}\label{conteq}
\nabla \cdot \boldsymbol{u}  = 0,
\end{align}
even in the ALE formulation.
 
Boundary and initial conditions are required to close the problem. In particular, the interface $\Gamma (t)$ is the surface of the moving cylinder. The quantities involved in the interface satisfy the following conditions
\begin{align}
& \boldsymbol{d}^g =y^C\boldsymbol{e}_y   &   \text{on} \quad & \Gamma (t) \\
& \boldsymbol{u} = \boldsymbol{u}^g =\dot{y}^C\boldsymbol{e}_y   &   \text{on} \quad & \Gamma (t) \\
& \boldsymbol{u}^g  = \boldsymbol{0} &   \text{on}  \quad & \partial \Omega \setminus \Gamma (t).
\end{align}
In the following, the coupling conditions at the interface are addressed.

\subsection{Coupling conditions at the interface}\label{couplingcond}
The coupling at the fluid-structure interface $\Gamma (t)$ is governed by,

\begin{itemize}
    \item Geometric compatibility: the domains must remain aligned at the interface,
    \item Kinematic continuity: the fluid velocity, mesh velocity, and structural velocity must be equal,
    \item Dynamic balance: the normal stresses of the fluid and the structure must be continuous.
\end{itemize}
The following equations summarize the above conditions.
{\RA{
\begin{align}
\boldsymbol{d}^g =y^C\boldsymbol{e}_y,  \quad \boldsymbol{u}^g =\dot{y}^C\boldsymbol{e}_y  ~~ \text{and} ~~  F_y = -\int_{\Gamma (t)}(\boldsymbol{\sigma}(\boldsymbol{x}, t)\cdot \boldsymbol{n})\cdot \boldsymbol{n}_y d\Gamma \qquad \text{on}\ \Gamma(t), 
\end{align}
}}
where,   the stress tensor for a Newtonian fluid is,
\begin{equation}
    \boldsymbol{\sigma}(\boldsymbol{x}, t) = -p(\boldsymbol{x}, t)\mathbf{I} + \mu \left(\nabla \cdot \boldsymbol{u}(\boldsymbol{x}, t) + \left( \nabla \cdot \boldsymbol{u}(\boldsymbol{x}, t) \right)^{\intercal} \right).
\end{equation}
Here, $\mathbf{I}$ is the identity tensor, $\bm{n}$  is the normal at the interface directed outward from the fluid domain.
{\RA{We finally point out that, in the framework of the two-way coupling implemented, the geometric and kinematic conditions involve information --- the wall displacements and velocities, respectively --- that goes from the solid problem to the fluid problem. Instead, the fluid force appearing in the dynamic condition represents the action of the fluid problem on the solid one.}}
\section{Numerical discretization of the full-order and reduced-order models}\label{NumDisFOMROM}
This section presents the discretization strategies adopted for both the full-order and reduced-order models, 
based on the Finite Volume Method (FVM) in the Arbitrary Lagrangian-Eulerian (ALE) framework.

\subsection{Numerical discretization of the full-order model}
\label{standardfv}
The aim of the FVM is to  discretize a system of partial differential equations written in integral form following \cite{moukalled2016finite}. In the present work, a 2-dimensional tessellation is used.  $N_{h}$ will represent the dimension of the full-order model (FOM) which is the number of cells of the discretized problem. 
In the following, the discretization methodology of the momentum and continuity equations will be addressed. In particular, the momentum and continuity equations will be solved using a segregated approach in the spirit of Rhie-Chow interpolation \cite{Rhie1983}. 

To approximate the problem by the use of the FVM, the domain $\Omega = \widehat{\Omega}(t) $ has to be subdivided through a tessellation $\mathcal{T}(t) = \{\Omega_i(t) \}_{i=1}^{N_h}$ so that every cell $\Omega_i(t)$ is a non-convex
polygon and $\bigcup_{i=1}^{N_h} \Omega_i(t) = \Omega (t)$ and $\Omega_i(t) \cap \Omega_j(t) = \emptyset \quad \forall i\neq j$. To simplify the notation, in the following,  $\Omega_i = \Omega_i(t)$ and $S_i = \partial \Omega_i(t)$. Where $S_i$ is the total surface related to cell $\Omega_i$.

We first consider the discretization of the continuity equation \ref{conteq} in the ALE framework. The finite volumes grid which moves in space must still obey the conservation law as pointed out in Tsui et al. \cite{tsui2013finite}, in which it is stated that as "the change in volume (area) of each control volume between time $t^n$ and $t^{n+1}$  must be equal to the volume (area) swept by the cell's boundary during $\Delta t = t^{n+1} -t^n$", namely
\begin{align}\label{gcl}
\frac{d }{d t}\int_{\Omega_i}d\Omega_i+  \int_{\partial \Omega_i}\boldsymbol{u}^g\cdot \boldsymbol{n} dS= 0.
\end{align}
For every control volume $\Omega_i = \Omega_i(t)$. Here, $\boldsymbol{n}$ is the outward unit normal vector on the boundary surface. The continuity equation written in integral form on control volume $\Omega_i$ reads

\begin{align}\label{eq:continuity_integral}
\frac{d }{d t}\int_{\Omega_i}\rho d\Omega_i+  \int_{\partial \Omega_i} \rho \left(\boldsymbol{u}- \boldsymbol{u}^g\right)\cdot \boldsymbol{n} dS= 0. 
\end{align}
Multiplying \cref{gcl} by the {\RB{density}} $\rho$ and making use of the incompressibility constraint yields,
\begin{align}
    \int_{\partial \Omega_i} \boldsymbol{u}^g \cdot \boldsymbol{n} dS = 0
\end{align}
which is the integral version of continuity equation \ref{conteq}, and shows that there is no need to consider the grid velocity in its discretization. 

In the framework of the {\RB {FVM}}, the momentum equation is written in its integral form for every cell of the tessellation
\begin{align}\label{IntegralMomentum}
    \int_{\Omega_i}\frac{\delta \boldsymbol{u}}{\delta t}d\Omega_i + \int_{\Omega_i}\nabla \cdot [\boldsymbol{u}\otimes(\boldsymbol{u} - \boldsymbol{u}_g)]d\Omega_i - \int_{\Omega_i}\nabla \cdot(\mu\nabla \boldsymbol{u})d\Omega_i  + \int_{\Omega_i}\nabla pd\Omega_i = 0.
\end{align}

In the FVM implementation here used, the volume integrals appearing in the momentum equation are computed by means of the so-called \emph{mid-point} quadrature rule. For such a reason, the unknowns of the algebraic problem resulting from the discretization are the values of the unknown fields $\boldsymbol{u}$ and $p$ at the centre of the cells. In particular, we denote with $\boldsymbol{u_i}$ and $p_i$ the velocity and pressure values at the centre of the generic cell $\Omega_i$. 

For a thorough discussion of the {\RB {FVM}}, the interested reader can refer to \cite{ferziger2002computational}. In the next sections,  the numerical treatment of each term in \cref{IntegralMomentum} is analysed and commented.

\subsubsection{The pressure gradient term} 
The pressure gradient term is discretized using Gauss's theorem.
\begin{equation}\label{pressure}
     \int_{\Omega_i}\nabla p d\Omega_i =  \int_{S_i}p d\boldsymbol{S}_i \approx \displaystyle\sum_{j \in S_i}\boldsymbol{S}_{ij}p_{ij}, 
\end{equation} 
where  $\boldsymbol{S}_{ij}$ is the oriented surface dividing the two neighbour cells $\Omega_i$ and $\Omega_j$, and $p_{ij}$ is the pressure evaluated at the centre of the face $\boldsymbol{S}_{ij}$. Also in the case of surface integrals, we make use of the mid-point quadrature rule.
\subsubsection{The convective term}
Making once again use of Gauss's theorem, the convective term is discretized as 
\begin{align}
    \int_{\Omega_i}\nabla \cdot [\boldsymbol{u}\otimes(\boldsymbol{u} - \boldsymbol{u}^g)]d\Omega_i  & = \int_{\Omega_i}\nabla \cdot (\boldsymbol{u}\otimes\boldsymbol{u})d\Omega_i - \int_{\Omega_i}\nabla \cdot (\boldsymbol{u}\otimes \boldsymbol{u}^g)d\Omega_i \\
    & = \int_{S_i}d\boldsymbol{S}_i\cdot (\boldsymbol{u} \otimes \boldsymbol{u}) -  \int_{S_i}d\boldsymbol{S}_i\cdot (\boldsymbol{u} \otimes  \boldsymbol{u}^g)\\ 
    &   = \displaystyle\sum_{j \in S_i}\boldsymbol{u}_{ij}F_{ij} - \displaystyle\sum_{j \in S_i}\boldsymbol{u}_{ij}({\boldsymbol{u}^g}_{ij} \cdot \boldsymbol{S}_{ij}) .
\end{align}
Here,  $\bm{u}_{ij}$ is the velocity evaluated at the centre of the face $\bm{S}_{ij}$, and $F_{ij} = \bm{u}_{ij}\cdot \bm{S}_{ij}$ is the flux of the velocity at the centre of the surface area vector $\bm{S}_{ij}$ sharing the cells $i$ and  $j$ oriented to the cell $j$.
This procedure underlines two considerations. The first one is that $\boldsymbol{u}_{ij}$ is not straightly available, all the variables of the numerical problem are evaluated at the centre of the cells. The values of the variables at the centre of the cell faces must then be obtained based on the values at the cell centres. There are of course several ways to obtain $\boldsymbol{u}_{ij}$ from $\boldsymbol{u}_{i}$. However, the basic idea behind them all is that the face value of each variable is obtained through interpolation of the cell centre values. The second clarification is about the fluxes. During an iterative process for the resolution of the equations, they are calculated using the velocity obtained at the previous step so that the non-linearity involved in such terms is quite naturally treated.
\subsubsection{The diffusion term }
The diffusion term is discretized as follows,
\begin{align}\label{diffusion}
    \int_{\Omega_i}\nabla \cdot  \left( \mu\nabla \boldsymbol{u} \right)d\Omega_i  = \mu_i\int_{S_i}d\boldsymbol{S}_i\cdot  \left( \nabla \boldsymbol{u} \right) \approx \displaystyle\sum_{j \in S_i} \mu_{ij} \boldsymbol{S}_{ij}\cdot (\nabla \boldsymbol{u})_{ij}, 
\end{align}
where $\mu_i$ is the viscosity of the $i-th$ cell, $\mu_{ij}$ is the viscosity evaluated at the centre of $\bm{S}_{ij}$, and $(\nabla \boldsymbol{u})_{ij}$ is the gradient of $\boldsymbol{u}_{ij}$ evaluated at the centre of $\bm{S}_{ij}$. 
As for the evaluation of the term $\boldsymbol{S}_{ij}\cdot (\nabla \boldsymbol{u})_{ij}$ in \cref{diffusion}, it involves the --- unknown --- gradient of the velocity at the face of the cell. Also in this case, we resort to interpolation to obtain such gradient. 
For orthogonal meshes,   the term  $\boldsymbol{S}_{ij}\cdot (\nabla \boldsymbol{u})_{ij}$ is computed as
\begin{align}\label{orthogonal-mesh}
\boldsymbol{S}_{ij}\cdot (\nabla \boldsymbol{u})_{ij} \approx 
\|\boldsymbol{S}_{ij}\|
\frac{\boldsymbol{u}_i- \boldsymbol{u}_j}{\|\boldsymbol{d}_{ij}\|}, 
\end{align}
where $\boldsymbol{d}_{ij}$ represents the vector connecting the centres of cells of index $i$ and $j$. We recall that a mesh is orthogonal if the line that connects two neighbouring cell centres is orthogonal to the face that divides the two cells. 

With non-orthogonal grids, the one-dimensional interpolation in \cref{orthogonal-mesh} loses accuracy, as the face centre is not lying on the line connecting the neighbouring cell centres. Thus, a correction term has to be added, and the interpolation corrected for non-orthogonality \cite{jasak1996error} reads
\begin{align}\label{non-orthogonal-mesh}
    \boldsymbol{S}_{ij}\cdot (\nabla \boldsymbol{u})_{ij} = \|\boldsymbol{\pi}_{ij}\|\frac{\boldsymbol{u}_i-\boldsymbol{u}_j}{\|\boldsymbol{d}\|} + \boldsymbol{k}_{ij} \cdot (\nabla \boldsymbol{u})_{ij}.
\end{align}
Herein, $\boldsymbol{S}_{ij} = \boldsymbol{\pi}_{ij} + \boldsymbol{k}_{ij}$ and  $\boldsymbol{\pi}_{ij}$ is chosen to be parallel to $\boldsymbol{S}_{ij}$ and $\boldsymbol{k}_{ij}$ to be orthogonal to $\boldsymbol{d}$. The term $(\nabla \boldsymbol{u})_{ij}$ is obtained through interpolation of the values of the gradient at the cell centres $(\nabla \boldsymbol{u})_i$ and $(\nabla \boldsymbol{u})_j$. 

\subsubsection{The PIMPLE algorithm}

Based on the aforementioned considerations, the discretized form of \cref{momentum} and  \cref{conteq} is written in the following matrix form
\begin{equation}\label{compactmatform}
    \left[ 
       \begin{array}{cc}
         \mathbf{A}_u   & \mathbf{B}_p  \\
           \nabla (\cdot) & \mathbf{0}
       \end{array}
       \right]
    \left[ 
       \begin{array}{c}
         \boldsymbol{u}_h \\
           \boldsymbol{p}_h
       \end{array}
\right]  =  \mathbf{0},
\end{equation}
where $ \boldsymbol{u}_h$ and $\boldsymbol{p}_h$ being the vectors where all $\boldsymbol{u}_i$  and  $p_i$ variables are collected respectively. With $\boldsymbol{u}_h \in \mathbb{R}^{dN_h}$ and $\boldsymbol{p}_h \in \mathbb{R}^{N_h}$.  ${N_h}$ is the number of cells in the mesh,  and  $d$ the  spacial dimension of the problem. Additionally, $\mathbf{A}_u$  is the matrix containing the terms related to velocity for the discretized momentum equation, $\mathbf{B}_p$ the matrix containing the terms related to pressure for the same equation, and $ \nabla (\cdot)$ the matrix representing the incompressibility constraint operator. The system matrix in \cref{compactmatform} has a saddle point structure which is usually difficult to invert using a coupled approach. For this reason, in this work we make use of a segregated approach, in which the momentum equation is solved with a tentative pressure and later corrected by exploiting the divergence-free constraint.

The segregated approach used in this work is based on the PIMPLE algorithm implemented in the OpenFOAM library \cite{jasak2007openfoam}. PIMPLE algorithm is a  mix of SIMPLE \cite{patankar1983calculation}, and  PISO \cite{issa1986solution} algorithms, and is mostly suitable for unsteady problems requiring high Courant number and/or dynamic mesh setups such as the one considered in this study.  

We here report a brief description of the PIMPLE algorithm. In particular, we mainly focus on the aspects that most affect the resolution of the online problem in the ROM strategy proposed in this work. 

As with the SIMPLE algorithm, PIMPLE is an iterative strategy that, at each time step, aims at converging to the correct pressure and velocity fields --- i.e. the ones respecting the continuity constraint. At each iteration, the first sub-step carried out is that of solving the discretized momentum equation, in which the pressure field is the one obtained at the previous iteration.  
In the following,  $\boldsymbol{u}_h^{n*}$ denotes the velocity which does not necessarily satisfy the continuity equation, while $\boldsymbol{u}_h^{n}$ does. 
The momentum matrix is divided into diagonal $\mathbf{A}$ and extra-diagonal parts $\mathbf{H}(\cdot)$
\begin{equation}\label{MomemtumMatrix}
  \mathbf{A}_u   \boldsymbol{u}_h^{n*} = \mathbf{A}\boldsymbol{u}_h^{n*} -  \mathbf{H}(\boldsymbol{u}_h^{n*}),
\end{equation}
with $n$ being an index to identify a generic iteration and $ \mathbf{A}_u$ satisfying the following relation 
\begin{align}\label{MomentumPredictor}
\mathbf{A}_u   \boldsymbol{u}_h^{n*} = -\mathbf{B}_p\boldsymbol{p}_h^{n-1}. 
\end{align}
By using \cref{MomemtumMatrix}, the momentum equation can be reshaped as follows
\begin{align}\label{velstar}
   \mathbf{A}\boldsymbol{u}_h^{n*} = \mathbf{H}(\boldsymbol{u}_h^{n*}) - \mathbf{B}_p\boldsymbol{p}_h^{n*} \Rightarrow & \boldsymbol{u}_h^{n*} =  \mathbf{A}^{-1}\mathbf{H}(\boldsymbol{u}_h^{n*})- \mathbf{A}^{-1}\mathbf{B}_p\boldsymbol{p}_h^{n-1}.
\end{align}
In a further iterative \textit{inner loop} small corrections to the velocity and pressure fields is introduced as
\begin{align}\label{corrections}
  \boldsymbol{u}_h^n = \boldsymbol{u}_h^{n*} + \boldsymbol{u}'  ~~~~~~~~~~ \boldsymbol{p}_h^n = \boldsymbol{p}_h^{n-1} + \boldsymbol{p}'.
\end{align}
The symbol $'$ denotes the corrections for both terms.
Inserting \cref{corrections} in \cref{velstar}, and rearranging terms yields
\begin{equation}\label{Mixterms}
   \boldsymbol{u}_h^n - \boldsymbol{u}' = \mathbf{A}^{-1}[\mathbf{H}(\boldsymbol{u}_h^{n}) - \mathbf{H}(\boldsymbol{u'})-\mathbf{B}_p\boldsymbol{p}_h^{n}+\mathbf{B}_p\boldsymbol{p}'].
\end{equation}
From \cref{Mixterms}, it is possible to deduce a relation between $\boldsymbol{u}'$ and $\boldsymbol{p}'$
\begin{equation}\label{uprime}
    \boldsymbol{u}' = \Tilde{\boldsymbol{u}}' - \mathbf{A}^{-1}\mathbf{B}_p\boldsymbol{p}', 
\end{equation}
with 
\begin{equation}\label{utildeprime}
    \Tilde{\boldsymbol{u}}' = \mathbf{A}^{-1}\mathbf{H}(\boldsymbol{u}').
\end{equation}
As the following relation holds thanks to \cref{MomentumPredictor} :
\begin{equation}
   \boldsymbol{u}_h^n = \mathbf{A}^{-1}[\mathbf{H}(\boldsymbol{u}_h^{n}) -\mathbf{B}_p\boldsymbol{p}_h^{n}].
\end{equation}
With the use of \cref{uprime} and the divergence operator $\nabla (\cdot)$ applied to $\boldsymbol{u}_h^{n}$ in \cref{corrections} knowing $\boldsymbol{u}'$ from \cref{uprime},  one obtain an equation that directly relates $\boldsymbol{p}'$ and $\boldsymbol{u}_h^{n*}$:
\begin{equation}
    [\nabla (\cdot)]\left(\mathbf{A}^{-1}\mathbf{B}_p\boldsymbol{p}'\right) = [\nabla (\cdot)]\boldsymbol{u}_h^{n*} + [\nabla (\cdot)]\Tilde{\boldsymbol{u}}'.
\end{equation}
Which is basically the discretized Poisson equation for pressure (PPE) expressed in terms of the velocity and pressure corrections.
In the SIMPLE algorithm, the velocity correction $\Tilde{\boldsymbol{u}}'$ is unknown as $\mathbf{H}(\boldsymbol{u}')$, hence neglected implying the following relation:  
\begin{equation}\label{relvelstarandpprime}
    [\nabla (\cdot)]\left(\mathbf{A}^{-1}\mathbf{B}_p\boldsymbol{p}'\right) = [\nabla (\cdot)]\boldsymbol{u}_h^{n*}.
\end{equation}
Therefore, $\boldsymbol{p}'$ is expressed as the only function of $\boldsymbol{u}_h^{n*}$ in \cref{relvelstarandpprime}. Then the corrected pressure is entered again in \cref{velstar} in order to obtain a new velocity field $\boldsymbol{u}_h^{n*}$ and repeat the procedure until the pressure correction falls below a given \textit{tolerance} and the velocity satisfy both the continuity and momentum equation. 

As the $\Tilde{\boldsymbol{u}}'$ is neglected, the SIMPLE algorithm converges slowly and is used mainly for steady-state simulations. Furthermore, to avoid instabilities, relaxation factor $\alpha_p$ and $\alpha_u$ are introduced in the computation of $\boldsymbol{p}_h^{n}$ and $\boldsymbol{u}_h^{n*}$ as follows:
\begin{align}
   &  \boldsymbol{p}_h^{n}  = \boldsymbol{p}_h^{n-1} + \alpha_p\boldsymbol{p}',\\
   & \boldsymbol{u}_h^{n*}  =  \mathbf{A}^{-1}\mathbf{H}(\boldsymbol{u}_h^{n*})- \alpha_u \mathbf{A}^{-1}\mathbf{B}_p\boldsymbol{p}_h^{n-1}.
\end{align}

The PISO algorithm comes to play to speed up the convergence after neglecting $\Tilde{\boldsymbol{u}}'$ and \textit{computing the pressure correction}  $\boldsymbol{p}'$  using \cref{uprime}. $\boldsymbol{u}'$ is  computed  as follows:
\begin{equation}\label{1stpresscorr}
    \boldsymbol{u}' = - \mathbf{A}^{-1}\mathbf{B}_p\boldsymbol{p}'.
\end{equation}
Allowing the computation of $\Tilde{\boldsymbol{u}}'$ using \cref{utildeprime}. One defines a second velocity correction equation mirroring \cref{uprime} as follows:
\begin{equation}\label{PISO}
    \boldsymbol{u}'' = \Tilde{\boldsymbol{u}}' - \mathbf{A}^{-1}\mathbf{B}_p\boldsymbol{p}''. 
\end{equation}
As $\boldsymbol{u}''$ in \cref{PISO} satisfy the continuity equation, one define  also a second pressure correction equation as:
\begin{equation}\label{2ndpresscorr}
      [\nabla (\cdot)]\left(\mathbf{A}^{-1}\mathbf{B}_p\boldsymbol{p}'' \right) =   [\nabla (\cdot)]\Tilde{\boldsymbol{u}}'.
\end{equation}
To sum up, what the PISO algorithm does more than the SIMPLE algorithm is to add a second inner loop to correct pressure and velocity. This speeds up the convergence, allowing this algorithm to be used in a transient simulation.  Following the procedure described by \cref{1stpresscorr,PISO,2ndpresscorr} further corrections steps can be  added, increasing both the algorithm's convergence and computational cost.

\subsection{Rigid body motion and mesh motion strategy}\label{motionstrategies}
In the resolution algorithm, the vertical displacement of the cylinder is accounted for by interfacing the fluid solver with Newton’s second law \cref{ode} for the rigid cylinder, written in the global inertial reference frame. Such a second-order ODE is in this work solved using the Symplectic 2nd-order explicit time-integrator for solid-body motion \cite{dullweber1997symplectic}.
After the linear displacement and velocity of the cylinder have been computed, the position of the cylinder boundary and the corresponding Dirichlet boundary datum for the velocity are updated in the fluid dynamic solver. Of course, due to rigid body kinematics consideration, the velocity prescribed on the cylinder boundary is uniform and equal to the velocity obtained solving \cref{ode}. To complete the interface between the structural --- rigid body --- solver and the fluid dynamic solver, the fluid grid must be updated to adjust to the boundary nodes motion while preserving good cells quality. 
 
In this work, the mesh deformation technique used is the so-called Slerp (Spherical Linear Interpolation) given its capability of deforming volumetric grids in presence of --- translational and rotational --- rigid boundary motions \cite{jasak2007openfoam}. For some alternatives to the mentioned mesh motion strategy, the interested reader may refer to \cite{MM2016, Jasak2007AutomaticMM, deBoer2007} and to references therein for more details. Finally, a finite difference scheme is used to obtain the grid velocity field $\boldsymbol{u}^g$ from the computed grid deformation field.

Finally, the essential steps of the PIMPLE algorithm combined  mesh motion strategy described in this section are reported in Algorithm~\ref{alg:pimple}.  In  Algorithm~\ref{alg:pimple},  the iterations within one time-steps are called \textit{outer iterations}, they are performed in an \textit{outer loop} in which the structural problem is solved, the cylinder position and the fluid mesh are updated, along with the coefficients and the source matrix of the discretized equations. The operations performed on linear systems with fixed coefficients are called instead \textit{inner iterations} and they occur in the so-called \textit{inner loop}.

\begin{algorithm}
\caption{PIMPLE algorithm with dynamic mesh.}
\label{alg:pimple}
\SetKwInOut{Input}{Input}
\SetKwInOut{Output}{Output}
\Input{\text{Initial fields $\boldsymbol{u}_h^{n*}$, $\boldsymbol{p}_h^{n-1}$, and ${\boldsymbol{d}^g}^0$}  \Comment{${\boldsymbol{d}^g}^0$ is the initial node displacement};} 
\Output{\text{$\boldsymbol{u}_h^{n}$, $\boldsymbol{p}_h^{n}$, and ${\boldsymbol{d}^g}^n$;}}
\While{ $t \leq t_{end}$ }
 {
      \While{$\text{ {\RB{Number of}} outer corrections $\geq$ 2}  ~~ \text{and} ~~\text{Tol} \geq \text{maxTol} $}
      {
          Compute the forces; \Comment{Using $\boldsymbol{u}_h^{n*}$, $\boldsymbol{p}_h^{n-1}$}\;
          Solve the rigid body problem \cref{ode}  \Comment{To obtain the new cylinder's position}\;
          Solve the mesh motion problem \Comment{To obtain ${\boldsymbol{d}^g}^n$ }\;
          $\mathbf{A}_u \boldsymbol{u}_h^{n*}$  \Comment{Assembling the momentum matrix \cref{MomemtumMatrix}}\;
          Solve $\mathbf{A}_u \boldsymbol{u}_h^{n*} = -\mathbf{B}_p\boldsymbol{p}_h^{n-1}$  \Comment{Momentum predictor \cref{MomentumPredictor} to obtain $\boldsymbol{u}_h^{n*}$}\;
          $ [\nabla (\cdot)]\left(\mathbf{A}^{-1}\mathbf{B}_p\boldsymbol{p}'\right) = [\nabla (\cdot)]\boldsymbol{u}_h^{n*}$ \Comment{Assembling the matrix of PPE \cref{relvelstarandpprime}}\;
          Solve $ [\nabla (\cdot)]\left(\mathbf{A}^{-1}\mathbf{B}_p\boldsymbol{p}'\right) = [\nabla (\cdot)]\boldsymbol{u}_h^{n*}$ \Comment{PPE to obtain $\boldsymbol{p}'$} \;
          $\boldsymbol{u}' \gets - \mathbf{A}^{-1}\mathbf{B}_p\boldsymbol{p}'$\Comment{ Momentum corrector \cref{1stpresscorr}} \;
          \While{$ \text{ {\RB {Number of}} inner corrections }$ } 
          {   
              $[\nabla (\cdot)]\left(\mathbf{A}^{-1}\mathbf{B}_p\boldsymbol{p}''\right) = [\nabla (\cdot)]\Tilde{\boldsymbol{u}}'$  \Comment{Assembling the matrix for PPE  
                \cref{2ndpresscorr}}\;
              Solve $ [\nabla (\cdot)]\left(\mathbf{A}^{-1}\mathbf{B}_p\boldsymbol{p}''\right) = [\nabla (\cdot)]\Tilde{\boldsymbol{u}}'$ \Comment{Recursively to obtain $\boldsymbol{p}''$}\;
              $\boldsymbol{u}' \gets \Tilde{\boldsymbol{u}}' - \mathbf{A}^{-1}\mathbf{B}_p\boldsymbol{p}''$: \Comment{Momentum corrector \cref{PISO}} \;
              }
              $\boldsymbol{u}_h^{n*} \gets \boldsymbol{u}'$ \;
              $\boldsymbol{p}_h^{n-1}\gets \boldsymbol{p}_h^{n-1} + \boldsymbol{p}'$ \;
       }
}
\end{algorithm}

\subsection{The reduced-order problem}\label{rompro}
In this section, we recall the details of the proper orthogonal decomposition, the concept of  POD-Galerkin projection,  and POD-Interpolation (POD-RBF) using radial basis networks.
\subsubsection{The Proper orthogonal decomposition}\label{POD}
The Proper Orthogonal Decomposition (POD) is used to construct the low-dimensional space. The POD is a compression technique where a set of numerical realizations (in time or parameter space) is reduced into a number of orthogonal basis (spatial modes) that capture the essential information suitably combined from previously acquired system data \cite{anttonen2005applications}. 

This work applies the POD to a group of realizations called snapshots. It consists of computing a certain number of full-order solutions $\boldsymbol{s}_i = \boldsymbol{s}(t_i)$ where $t_i \in \boldsymbol{T}$ for $i = 1, \cdots, N$. $\boldsymbol{T}$ being the training collection of a certain number $N$ of the time values, to obtain a maximum amount of information from this costly stage to be employed later on for a cheaper resolution of the problem. Those snapshots can be assembled at the end of the resolution into a so-called \emph{snapshot matrix} $\boldsymbol{S} \in \mathbb{R}^{N_h\times N}$ defined as
\begin{align}
    \boldsymbol{S} & = \left[\boldsymbol{s} (\boldsymbol{x},t_1), \dots, \boldsymbol{s} (\boldsymbol{x},t_{N}) \right].
\end{align}
The idea is to compute the ROM solution that can minimize the error denoted here by $E^{ROM}$ see \cref{Erom} between
the obtained realization of the problem and its high-fidelity counterpart. 
In the POD-Galerkin scheme, the reduced-order solution is represented as 
\begin{align}\label{seqq}
\boldsymbol{s} (\boldsymbol{x},t) \approx  \boldsymbol{s}^{ROM} (\boldsymbol{x},t) = \displaystyle\sum_{i=1}^{N_r}a_i(t)\boldsymbol{\phi}_i(\boldsymbol{x}).
\end{align}
Where $N_r \ll N_h$ ($N_h$ is the number of cells in the computational domain) is a predefined number,  $\boldsymbol{\phi}_i$ is a generic pre-calculated orthonormal function depending only on the space while $a_i(t)$ is the temporal modal coefficients  satisfying the conditions  

\begin{align}
{\RB {
    a_j(t)  =   \left(\boldsymbol{\phi}_j,  \boldsymbol{s} (\boldsymbol{x},t) \right)_{L^2 (\Omega)}, ~~
    \boldsymbol{\phi}_j^T \bm{M} \boldsymbol{\phi}_i = \bm{\delta}_{ij}.
    }
}
\end{align}
$\bm{M}$ being  the mass matrix defined by the chosen inner product. In the case of $L_2$-norm and FVM $\bm{M}$ is a diagonal matrix containing the cell volumes. 
The best performing functions $\boldsymbol{\phi}_i$ in this case, are  the ones minimizing the  $L^2$-norm error $E$ between all the reduced-order solutions $\boldsymbol{s}^{ROM}_i$, $i = 1, \cdots, N$ and their high fidelity counterparts, namely
\begin{align}\label{Erom}
    E = \displaystyle\sum_{i=1}^{N}\|\boldsymbol{s}^{ROM}_i - \boldsymbol{s}_i\|_{L^2(\Omega)} =  \displaystyle\sum_{i=1}^{N}\|\boldsymbol{s}_i - \displaystyle\sum_{i=1}^{N_r}\left(\boldsymbol{s}_i, \boldsymbol{\phi}_i\right)_{L^2(\Omega)}\boldsymbol{\phi}_i \|_{L^2(\Omega (t_0))}.
\end{align}
$\Omega (t_0)$ being the reference configuration of the computational domain in the case of grid motion. Note that the projection is performed with respect to $L_2(\Omega(t))$ while POD is computed with respect  to $L_2(\Omega(t_0))$. 
It can be shown \cite{Kunisch2002} that solving a minimization problem based on \cref{Erom} is equivalent to solving the eigenvalue problem 
\begin{align}\label{eigenproblem}
    \boldsymbol{C} \mathbf{V} = \mathbf{V}\boldsymbol{\lambda}.
\end{align}
 $\boldsymbol{C} \in \mathbb{R}^{N\times N}$ being the correlation matrix between all the different training solutions of the snapshot matrix $\boldsymbol{S}$, $\mathbf{V} \in \mathbb{R}^{N\times N}$ is the matrix whose columns are the eigenvectors, and $\boldsymbol{\lambda} \in \mathbb{R}^{N\times N}$  is a diagonal matrix whose diagonal entries are the eigenvalues.  The entries of the correlation matrix are defined as follows
\begin{align}\label{correlation}
\boldsymbol{C}_{ij} = \left(\boldsymbol{s}_i, \boldsymbol{s}_j \right)_{L^2\Omega (t_0))},
\end{align}
using a POD strategy, the required basis functions are obtained through the resolution of the eigenproblem mentioned in \cref{eigenproblem}, obtained with the method of snapshots by solving \cref{Erom}. The required basis functions  are then computed from the eigenvalues and eigenvectors in \cref{eigenproblem} as 
\begin{eqnarray}
    \boldsymbol{\phi}_i = \displaystyle\dfrac{1}{N\sqrt{\lambda_i}}\displaystyle\sum_{j=1}^{N}\boldsymbol{s}_jV_{ji} ~~~~~ \forall i = 1, \cdots, N.
\end{eqnarray}
All the basis functions are collected into a single matrix:
\begin{equation}
      \boldsymbol{\Phi} = \left[\boldsymbol{\phi}_1, \cdots, \boldsymbol{\phi}_{N_r} \right] \in \mathbb{R}^{N_h \times N_r}.
\end{equation}
Which is used to project the high fidelity problem onto the reduced subspace so that the final system dimension is $N_r$. 
In the framework of ROM analysis, the \emph{offline phase} consists in carrying out multiple resolutions of the FOM problem, collecting the snapshots to assemble the snapshot matrix, computing the modes through eigenproblem \cref{eigenproblem}, and projecting the FOM problem onto the reduces subspace. Once this computationally expensive phase is done, the procedure results in a solution system (the \emph{online phase}) characterized by a small amount of unknowns, and by computational cost that is much lower than the original problem. In the next sections, we will provide details on the way we adapted the Galerkin projection procedure to the present fluid-structure interaction application.

\subsubsection{Reduced-PIMPLE algorithm for incompressible laminar flows}\label{sectredpimple}
This part assumed that the  high fidelity  model is discretized and written in the  form    
\begin{align}\label{discritizefullmodel}
\mathbf{A}_u \boldsymbol{u}_h = \boldsymbol{b}_u,  ~~~~~~	\mathbf{B}_p \boldsymbol{p}_h = \boldsymbol{b}_p.
\end{align}
With $\mathbf{A}_u \in \mathbb{R}^{dN_h\times dN_h}$, $\boldsymbol{u}_h \in \mathbb{R}^{dN_h}$, $\mathbf{B}_p \in \mathbb{R}^{N_h\times N_h}$, $\boldsymbol{p}_h \in \mathbb{R}^{N_h}$, $\boldsymbol{b}_u \in \mathbb{R}^{N_h}$,  and $\boldsymbol{b}_p \in \mathbb{R}^{N_h}$ are  defined  in \cref{compactmatform}.  
As already mentioned, $N_h$ is the number of control volumes (cells) in the mesh, and $d=2$ is the space dimension. 
The following equation introduces the resulting reduced model expansions of the velocity and pressure fields respectively
$\boldsymbol{u}_h (\boldsymbol{x},t) \approx \boldsymbol{u}_r (\boldsymbol{x},t)$ and $\boldsymbol{p}_h(\boldsymbol{x},t) \approx p_r(\boldsymbol{x},t)$, with
\begin{align}\label{ueqq}
\boldsymbol{u}_r (\boldsymbol{x},t) =  \displaystyle\sum_{i=1}^{N_u}a_i(t)\boldsymbol{\phi}_i(\boldsymbol{x}) = \boldsymbol{\Phi} \boldsymbol{a}^T,~~~~~ p_r(\boldsymbol{x},t) = \displaystyle\sum_{i=1}^{N_p}b_i(t)\boldsymbol{\xi}_i(\boldsymbol{x}) = \boldsymbol{\Xi} \boldsymbol{b}^T.
\end{align}
Herein,  $a_i(t)$,  and  $b_i(t)$ are temporal modal coefficients;  $\boldsymbol{\phi}_i$ and $\boldsymbol{\xi}_i$ are the basis functions of POD modes of the velocity and pressure fields stored  respectively in $\boldsymbol{\Phi} \in \mathbb{R}^{dN_h\times N_u}$ and $\boldsymbol{\Xi} \in \mathbb{R}^{N_h\times N_p}$ with $N_u$ and $N_p$ being the numbers of basis functions selected for the prediction of the velocity and pressure solutions respectively.  $\boldsymbol{a} \in \mathbb{R}^{N_u}$ and $\boldsymbol{b} \in \mathbb{R}^{N_p}$ are  the vectors containing the temporal coefficients for the velocity and pressure, respectively.  
We point out that for the construction of the reduced basis spaces, the POD strategy discussed in \cref{rompro} is used on the snapshot matrices of the velocity and pressure fields separately, in order to obtain two different families of reduced basis functions.
\begin{align}\label{basismatrices}
\boldsymbol{\Phi}= \left[\boldsymbol{\phi}_1,\dots, \boldsymbol{\phi}_{N_u} \right] \in \mathbb{R}^{dN_h\times N_u}, ~~~~~~
\boldsymbol{\Xi} = \left[\boldsymbol{\xi}_1, \dots, \boldsymbol{\xi}_{N_p} \right] \in \mathbb{R}^{N_h\times N_p}.
\end{align}
The linear systems in \cref{discritizefullmodel} are projected respectively into the low-dimensional using in \cref{basismatrices}. Thus, the following relation holds: 
\begin{align}\label{redsystems}
 \boldsymbol{A}_u^r\boldsymbol{a} = \boldsymbol{b}_u^r, ~~~~~~~~~~~~~~~~  \boldsymbol{A}_p^r\boldsymbol{b} = \boldsymbol{b}_p^r.
\end{align}
Where $\boldsymbol{A}_u^r = \boldsymbol{\Phi}^T\mathbf{A}_u\boldsymbol{\Phi} \in \mathbb{R}^{N_r^u\times N_r^u}$,    $\boldsymbol{A}_p^r = \boldsymbol{\Xi}^T\mathbf{A}_p\boldsymbol{\Xi} \in \mathbb{R}^{N_r^p\times N_r^p}$, $\boldsymbol{b}_u^r = \boldsymbol{\Phi}^T\boldsymbol{b}_u \in \mathbb{R}^{N_r^u}$, and $\boldsymbol{b}_p^r = \boldsymbol{\Xi}^T\boldsymbol{b}_p \in \mathbb{R}^{N_r^p}$.
The resulting systems in \cref{redsystems} can be solved using any method for dense matrices. In this work, we used the Householder rank-revealing QR decomposition of a matrix with full pivoting implemented in the Eigen library \cite{eigenweb}. The whole idea here is to rely on a method capable of being as coherent as possible with respect to the high-fidelity Algorithm~\ref{alg:pimple} discussed earlier. The main steps of the reduced algorithm
for incompressible laminar flows are summarized in Algorithm~\ref{alg:redpimple}. 
We point out that in the current version of the algorithm, steps 7 and 10 have to be carried out on the full order grid, without a consequent impact on the reduced order model computational cost. 
In \cref{reduced-pimple-steps}, a step-by-step flowchart to  clarify Algorithm~\ref{alg:redpimple} is reported.

\begin{algorithm}
\caption{Reduced-PIMPLE algorithm with dynamic mesh}
\label{alg:redpimple}
\SetKwInOut{Input}{Input}
\SetKwInOut{Output}{Output}
\Input{\text{$\boldsymbol{u}_h^{n*}$, $\boldsymbol{p}_h^{n-1}$, ${\boldsymbol{d}^g}^0$, $\boldsymbol{\Phi}$, $\bm{\Psi}$, and $\boldsymbol{\Xi}$;}}
\While{ $t \leq t_{end}$ }
 {
      \While{$\text{ {\RB {Number of}} outer corrections $\geq$ 2}  ~~ \text{and} ~~\text{Tol} \geq \text{maxTol} $}
      {
          Compute the forces; \Comment{Using $\boldsymbol{u}_h^{n*}$, $\boldsymbol{p}_h^{n-1}$}\;
          Solve the rigid body problem \cref{ode}  \Comment{To obtain the new cylinder's position $y^C_{new}$}\;
          Compute $\boldsymbol{c} = RBF (y^C_{new})$ \cref{rbfnns}\;    
          Reconstruct ${\boldsymbol{d}^g}^n  = \bm{\Psi} \boldsymbol{c}^T$ \cref{POD_displacements}\;
          $\mathbf{A}_u \boldsymbol{u}_h^{n*}  = \mathbf{b}_u$ \Comment{Assembling the momentum matrix \cref{MomentumPredictor}}\;
          Solve $\boldsymbol{\Phi}^T \mathbf{A}_u \boldsymbol{\Phi} \boldsymbol{a}^{*} =\boldsymbol{\Phi}^T \mathbf{b}_u$ \Comment{To obtain $\boldsymbol{a}^{*}$ \text{with} $\mathbf{b}_u =-\mathbf{B}_p\boldsymbol{p}_h^{n-1}$} \;
          Reconstruct $\boldsymbol{u}_h^{n*}$\  \Comment{Using $\boldsymbol{a}^{*}$}\;
          $ [\nabla (\cdot)]\left(\mathbf{A}^{-1}\mathbf{B}_p\boldsymbol{p}'\right) = [\nabla (\cdot)]\boldsymbol{u}_h^{n*}$\Comment{Assembling the matrix of PPE \cref{relvelstarandpprime}}\;
          Solve $\boldsymbol{\Xi}^T \mathbf{A}_p \boldsymbol{\Xi} \boldsymbol{b}' = \boldsymbol{\Xi}^T\boldsymbol{b}_p$ \Comment{To obtain $\boldsymbol{b}'$}\;
          Reconstruct $\boldsymbol{p}'$\  \Comment{Using $\boldsymbol{b}'$}\;
          $\boldsymbol{u}' \gets - \mathbf{A}^{-1}\mathbf{B}_p\boldsymbol{p}'$ \Comment{Momentum corrector \cref{1stpresscorr}} \;
          \While{\text{ {\RB {Number of}} inner corrections } } 
          {   
              $[\nabla (\cdot)]\left(\mathbf{A}^{-1}\mathbf{B}_p\boldsymbol{p}''\right) = [\nabla (\cdot)]\Tilde{\boldsymbol{u}}'$  \Comment{Assembling the matrix for PPE  
                \cref{2ndpresscorr}}\;
              Solve $\boldsymbol{\Xi}^T \mathbf{A}_p \boldsymbol{\Xi} \boldsymbol{b}'' = \boldsymbol{\Xi}^T\boldsymbol{b}_p$ \  \Comment{Recursively to obtain $\boldsymbol{b}''$ \text{where} $\boldsymbol{b}_p = [\nabla (\cdot)]\Tilde{\boldsymbol{u}}'$}\;
              Reconstruct $\boldsymbol{p}''$\  \Comment{Using $\boldsymbol{b}''$}\;
              $\boldsymbol{u}' \gets \Tilde{\boldsymbol{u}}' - \mathbf{A}^{-1}\mathbf{B}_p\boldsymbol{p}''$: \Comment{Momentum corrector \cref{PISO}} \;
           }
           $\boldsymbol{u}_h^{n*} \gets \boldsymbol{u}'$ \;
           $\boldsymbol{p}_h^{n-1}\gets \boldsymbol{p}_h^{n-1} + \boldsymbol{p}'$ \;
       }
}
\Output{\text{ $\boldsymbol{u}_h^{n}$, $\boldsymbol{p}_h^{n}$, and ${\boldsymbol{d}^g}^n$;}}
\end{algorithm}

\subsubsection{POD with interpolation for mesh motion prediction}
\label{pod-i}
This section presents a method to reduce the computational cost associated with the mesh motion part in the system. {\RA{On one hand, it must be first remarked that the structural problem, at the numerical level, is characterized only by a single degree of freedom. For such a reason, the same --- inexpensive --- structural equations used in the full order model are solved also at the reduced level. On the other hand, the computational cost associated with the mesh morphing strategy described in Section \ref{motionstrategies} is considerably higher and is here reduced by means of a POD-based algorithm. }}Along with a reduction of resolution system degrees of freedom, the advantage of this methodology is to make the online equations independent of the specific equations solved at the full order levels to compute the mesh deformation. As will be discussed,  the methodology combines proper orthogonal decomposition with radial basis functions (RBF) networks applied to the grid nodes displacement field.  

So, the first step of the mesh deformation reduction strategy is that of computing the POD modes of the grid nodes displacement field. To this end, we assemble a snapshot matrix with the grid nodes displacements obtained at different time steps

\begin{align}
    \boldsymbol{S}^g & = \left[\boldsymbol{d}^g (\boldsymbol{x},t_1), \dots, \boldsymbol{d}^g (\boldsymbol{x},t_{N}) \right].
\end{align}
As in the case of pressure and velocity unknowns, the matrix $\boldsymbol{S}^g$ is then processed to obtain a correlation matrix using \cref{correlation} and, from the solution of an eigenvalue problem as in \cref{eigenproblem}, a set of POD modes. The reduced-order solution for the {\RB {grid nodes’ displacement field}} is then represented as, 
\begin{align}\label{POD_displacements}
\boldsymbol{d}^g (\boldsymbol{x},t) \approx   \displaystyle\sum_{i=1}^{\RB {N_{\text{pD}}}}c_i(t)\boldsymbol{\psi}_i(\boldsymbol{x}),
\end{align}
where ${\RB {N_{pD}}}$ is the amount of modes considered for the grid displacement field. Along with the modal functions $\boldsymbol{\psi}_i(\boldsymbol{x})$, the solution of eigenvalue \cref{eigenproblem} provides the values of the modal coefficients in correspondence with each time step included in the snapshot matrix. For such a reason, a natural choice for computing the grid deformation at time instants not included in the snapshots would be that of interpolating the modal coefficients based on the time variable. However, given the fact that the grid nodes displacement is induced by the --- rigid --- translation of the cylinder boundary $y^C$, a more meaningful way to obtain the modal coefficients at each time step is to consider that
\begin{equation}
c_i(t)=\widehat{c}_i(y^C(t))\quad i=1,\dots, {\RB {N_{\text{pD}}}},    
\end{equation}
and interpolate the $c_i$ values based on the cylinder vertical displacement variable --- as obtained at each of the time steps at which the solution snapshots have been collected.

In this work, the interpolation step of the data driven POD strategy used for the reduction of the grid nodes displacement field is carried out by means of the Radial Basis Function \cite{LAZZARO2002521} method. In the present framework, modal coefficient $c_i$ at a generic value $y^C$ is obtained evaluating the expression
\begin{equation}\label{rbfnns}
    c_i(y^C)  = \displaystyle\sum_{k=1}^{N}w_k{\RB {\varphi}}(|| y^C - y^C_k ||), 
\end{equation}
in which ${\RB {\varphi}}:\mathbb{R}\rightarrow\mathbb{R}$ is in the so-called \emph{radial basis}.  ${\RB {\varphi}}$ is a function of the Euclidean distance. In the present case it is a one dimensional function, but in more general cases it maps the $m$ dimensional parameter space in $\mathbb{R}$. The weights $w_k$ appearing in \cref{rbfnns} are determined imposing the interpolation condition at the snapshots, in which the modal coefficients are known from \cref{eigenproblem}. The conditions used are then

\begin{equation}\label{rbfnns_colloc}
    (c_i)_j=c_i(y^C_j)  = \displaystyle\sum_{k=1}^{N}w_k{\RB {\varphi}}(|| y_j^C - y^C_k ||) \quad j=1,\dots,N 
\end{equation}
resulting in the system
 \begin{equation}\label{rbfweights}
   \boldsymbol{c}   = \mathbf{G} \boldsymbol{w}^T, 
 \end{equation}
in which $ \mathbf{G} = (g_{kj}) = {\RB {\varphi}}(|| y_j^C - y^C_k ||)$ is the Gram matrix. Once the weights have been computed solving system \cref{rbfweights} in the offline phase, in the online computations the modal coefficients are obtained evaluating \cref{rbfnns}.

We point out that in the present case the $c_i$ coefficients only depend on the scalar variable $y^C(t)$. In the case of multidimensional dependence, RBF interpolation can still be used with no algorithmic modifications. Thus, more complex cases in which the cylinder exhibits rigid motions with more translational and rotational degrees of freedoms, could still be treated with the methodology described. Even in the case in which the FSI problem involves a deformable body which alters the shape of one or more boundaries of the fluid domain, RBF 
could be still used to interpolate the fluid mesh nodal displacements based on the structural displacements reduced coefficients. 

Finally, we remark that the RBF interpolation can also be interpreted as a network in which $N$ is the number of neurons in the hidden layer, $y^C_k$ is the centre vector for neuron $k$, $\boldsymbol{w}_k$ being the weight of neuron $k$ in the linear output neuron. Given this analogy, we point out that different and more efficient networks can substitute RBF --- which has a $O(N^2)$ computational cost --- in future works.
\section{Numerical tests}
\label{resdiscuss} 
\begin{figure}[h!]
    \centering
    \begin{tabular}{cc}
          
  \ifpdf
  \resizebox{0.3\textwidth}{!}{
    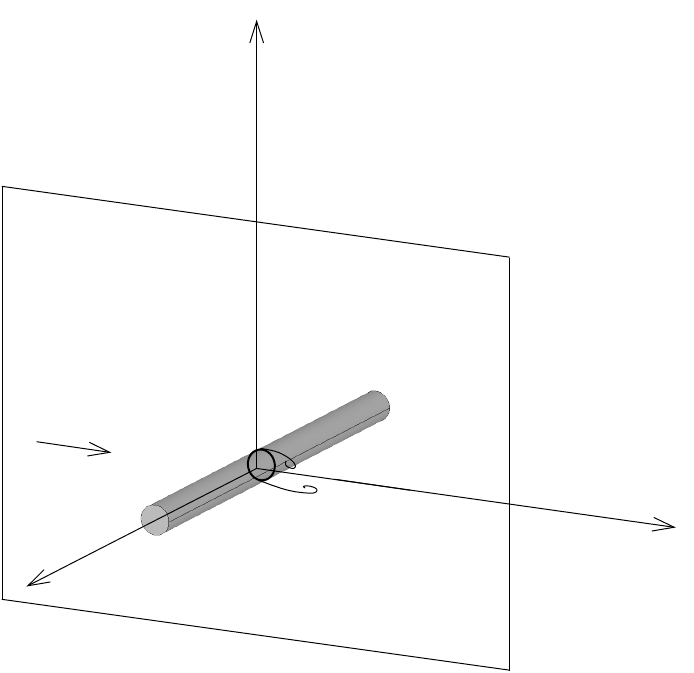
  }
  \else
  \resizebox{0.3\textwidth}{!}{
    \input{./images/cylinder_3D_2D.pstex_t}
  }
  \fi

         &  \includegraphics[width=0.4\textwidth]{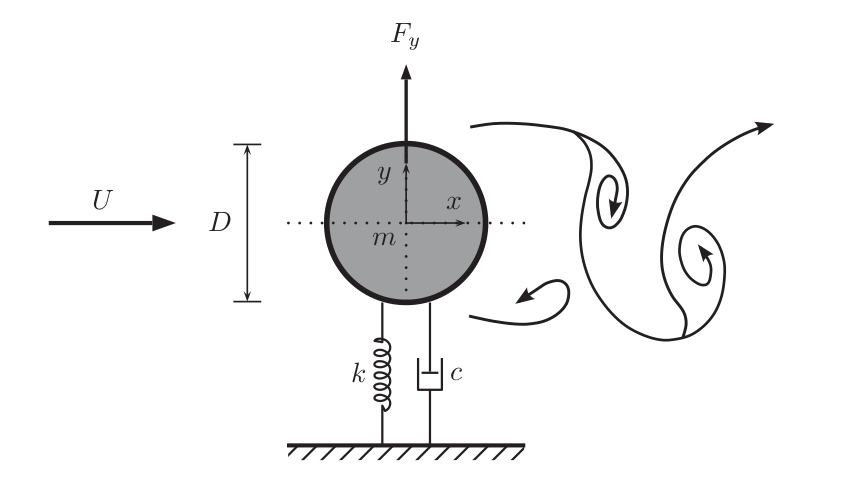}
         \\
         (a) & (b)
    \end{tabular}
    \caption{Problem setup: (a) Conceptual 3D schematic illustrating flow past a cylinder and vortex shedding. This is for visualization purposes only; the actual simulations are two-dimensional;
    (b) 2D computational domain used for the simulation. The mesh and boundary conditions are defined in this plane.}
    \label{fig:cross-flow-viv}
\end{figure}

As previously discussed, the physical problem considered in this work is that of an elastically mounted cylinder restrained to move in the transverse direction, as shown in Fig. \ref{fig:cross-flow-viv} (b). {\RA{It is important to stress that vortex shedding on elongated cylinders is an inherently two dimensional flow. In such conditions tip effects are in fact not affecting the flow around the vast majority of the body. As such, 2D simulations on sections of the cylinder located away from its tips are a good approximation of the local 3D flow field. This is illustrated in Fig. \ref{fig:cross-flow-viv}, where in picture (a) we show a three dimensional sketch of the cylinder, including a sectional plane normal to the cylinder axis and parallel to the inflow velocity $U_{\infty}$. Such a sectional plane, is the domain $\Omega$ in which --- as shown in Fig. \ref{fig:cross-flow-viv} (b) --- we solve a two dimensional fluid-structure interaction problem in which the cylinder circular section is assumed to be rigid. In this framework, the elastically mounted circle is also restrained to move in the transverse direction.}}

\subsection{Description of the configuration and boundary conditions}
The computational domain has a length of 34 $D$ and a width of 10$D$, where $D = \SI{1.0}{\metre}$ is the cylinder diameter. The cylinder is located at a 5$D$ distance from the inlet.   Fig. \ref{fig:mylabel} presents a view of the two-dimensional computational grid, both in its reference/initial  configuration, and in a deformed state caused by the $~0.4D$ vertical displacement of the cylinder. The grid features 11 644 cells (control volumes) and 24 440 points. 
The flow velocity at the inlet is $\boldsymbol{U}_{\infty} = (U_{in}, 0)$ with $U_{in} = \SI{1.0}{\metre\per \s}$, and the  physical viscosity $\nu = \SI{0.005}{\kilogram\per{\metre\s}}$. 
This corresponds to a Reynolds number of 200. Next, the boundary conditions are summarized in Table \ref{table:boundaryconditions}.
\begin{longtable}{ccccc}
&  Inlet & Sides & Outlet & Cylinder\\ 
 \hline
$\boldsymbol{u}$ & $\boldsymbol{u} = (1, 0)$ & $\boldsymbol{u} \cdot \boldsymbol{n} = 0$  & $\nabla \boldsymbol{u} \cdot \boldsymbol{n} = 0$ & $\boldsymbol{u} = (0, \dot{y}^C)^*$\\  
$p$ & $\nabla p\cdot \boldsymbol{n} = 0$ & $\nabla p\cdot \boldsymbol{n} = 0$  & {\RB {$p=0$}} & {\RB {$p=0$}}\\
$\boldsymbol{d}^g$ & $\boldsymbol{d}^g = \boldsymbol{0}$ &  $\boldsymbol{d}^g = \boldsymbol{0}$ & $\boldsymbol{d}^g = \boldsymbol{0}$ &  $\boldsymbol{d}^g=(0, y^C)^*$\\  
 \hline
 \caption{A summary of the boundary conditions imposed in the ALE fluid dynamic problem. Note that the $^*$ subscript indicates quantities that are computed by the rigid body structural solver.}
\label{table:boundaryconditions}
\end{longtable}
At the inlet boundary, non-homogeneous Dirichlet conditions are imposed on the velocity field, while a zero-gradient condition is applied to the pressure field. At the outlet, zero-gradient conditions are prescribed for the velocity field, and homogeneous Dirichlet conditions are imposed on the pressure. Along the lateral boundaries (top and bottom), zero-gradient conditions are applied to both velocity and pressure fields. On the cylinder surface, interface coupling conditions from the structural solver, as described in \cref{couplingcond}, are enforced.
\subsubsection{Linear solvers for the fluid}
The simulations are carried out using the PIMPLE  Algorithm~\ref{alg:pimple}. 
{\RA {Although the implicit time-stepping scheme (implicit Euler) used for computing the time derivative of the velocity field is unconditionally stable, we monitor the Courant–Friedrichs–Lewy (CFL) number to ensure adequate temporal resolution and solution accuracy, particularly in regions of rapid flow variation. 
In this simulation, the PIMPLE algorithm adapts the time step to maintain the maximum CFL number below a prescribed threshold of 0.5.}}For the spatial gradients, {\RB {a Gauss linear}} scheme has been employed. 
The convective term has been approximated with the Upwind scheme.  
Gauss linear scheme is used to approximate the diffusive term. 
The values of the relaxation factors $\alpha_u$, and $\alpha_p$ have been fixed at 0.7 and 0.3, respectively. 
One non-orthogonal corrector iteration is used to deal with the mesh's non-orthogonality. 
In addition, two pressure corrector and one momentum correctors are used in the simulations.
As for the linear solvers, a smoother Gauss-Seidel has been used for solving the momentum equation, and  GAMG (geometric-algebraic multi-grid) for solving the pressure equation. 

\subsubsection{Structural solver}
As mentioned, the structural model is represented by the second-order differential  \cref{ode} for rigid-body motion, here  solved using the Symplectic 2nd-order explicit time-integrator \cite{dullweber1997symplectic}. The mass of the cylinder considered in the numerical tests is $m=\SI{0.05}{\kilogram}$, the spring stiffness is $m=\SI{6.76e-2}{\newton \per \metre}$, which results in the natural frequency $f_n=\SI{0.185}{\hertz}$.
The cylinder to ground connection damping coefficient is $c=\SI{0.0135}{\kilogram \per \second}$. 
The flow and structure parameters are summarized in  Table~\ref{table:params}.  
\begin{longtable}{ccccc}
 $Re$ & $f_n$ [Hz] & $c$ [kg/s] & $k$[N/m] & $m$ [kg]\\ 
 \hline
 200 & 0.185 & 0.0135 & 6.76e-2 & 0.05\\  
 \caption{Simulation parameters}
\label{table:params}
\end{longtable}

\subsection{Performance evaluation and parametric analysis of the reduced-order model}
The primary objective of the present numerical study is to assess the ability of the ROM to accurately predict the flow fields associated with the final periodic regime. 
To this end, the ROM's solver's performance on a benchemark case involving laminar flow past a circular cylinder undergoing vortex-induced vibrations at a Reynolds number $Re=200$.
The high-fidelity data used to collect the snapshot set were obtained using the FVM implemented in the C++ open-source library  OpenFOAM\textsuperscript{\textregistered} \cite{ofsite,jasak2007openfoam}. 
At the reduced-order level, both the modal reduction and the assembly and solution of the reduced-order systems were performed  performed using the C++ open-source library ITHACA-FV (In real Time Highly Advanced Computational Applications for Finite Volumes) \cite{stabile2018finite, StabileHijaziMolaLorenziRozza2017}. 
ITHACA-FV is specifically designed to interface with FVM solvers available in OpenFOAM, enabling seamless integration and facilating the application of ROM techniques to  real-world industrial problems, where OpenFOAM is commonly employed.
Additionally, the SPLINTER C++ library \cite{SPLINTER} has been used to construct the Radial Basis Function (RBF) networks employed in this work.
For the present cross-flow cylinder test case, the FOM simulation was first run until a periodic regime was fully established. Once this steady periodic behavior was reached, the simulation was continued for an additional $30\ $s using a constant time step \SI{0.003}{\s}, during which the solution fields were exported at intervals of \SI{0.1}{\s}.

\subsubsection{Computational cost}
Table~\ref{tab:OffOnline} reports a comparison analysis of the full-order and reduced-order models execution times as the number of modes for the prediction of velocity, pressure, and grid nodes displacement fields are varied. This allows for evaluating the effect of the number of modes variation on the computational cost of the online phase.
\begin{longtable}[h!]{|c|c|c|}
	\hline
	\textbf{Stages} & \textit{\# of modes} & \textit{Time} [s]\\
 \hline
	\textbf{Offline PDE solution} &  - & 1745.66 \\ 
 \hline
    \multirow{3}{*}{\textbf{Online PDE solution}} & $N_u= 25, N_p=20$, $N_{\textit{pD}}=3$ & 672.65\\ 
       \cline{2-3}
       &    $N_u=N_p=20$,$N_{pD}=3$ & 660.882 \\ 
       \cline{2-3}
       & $N_u=20, N_p=15$, $N_{\textit{pD}}=3$ & 631.8676563\\ 
       \cline{2-3}
       &$N_u =20,N_p =10$ $N_{\textit{pD}}=3$ &  616.8310596\\ 
       \hline
       \caption{Offline and Online times comparison varying the number of modes}
\label{tab:OffOnline}
\end{longtable}

The offline stage comprises four steps: 
\begin{itemize}
    \item Computing the snapshot set via numerical approximation of the original high-dimensional system,
    \item Deriving the POD basis,
    \item Projecting the dynamics onto the low-rank subspace,
    \item Evaluating the  radial basis network's.
\end{itemize}
However, only the computational cost of the first step is reported in Table~\ref{tab:OffOnline} as it is the most expensive. 

The online cost corresponds to the time required to solve the surrogate model. 
The computational times reported in Table~\ref{tab:OffOnline}, indicate that the ROM solution achieves only a modest speed-up compared to the FOM solver. This limited improvement arises because, in the deforming domain applications like FSI simulations, the ROM systems matrices must be recomputed at each step using both velocity and pressure projection matrices as done in \cref{redsystems}. It is worth noting that these projection matrices depend on $N_h$, the full-order grid size. This process significantly offsets the theoretical gains from reducing the number of unknows and currently represents a major bottleneck for computational efficiency. Ongoing efforts to address this issue focus on hyper-reduction.
Despite the challenges, the main goal of the present work is to assess the accuracy of the proposed ROM framework.  In particular, it is important to verify that the hybrid approach - merging physics-based model reduction with data-driven approximations of flow and displacement fields - achieves sufficient accuracy.

\subsubsection{Reconstruction error}

\begin{figure}[h!]
\centering
\includegraphics[width=\textwidth]{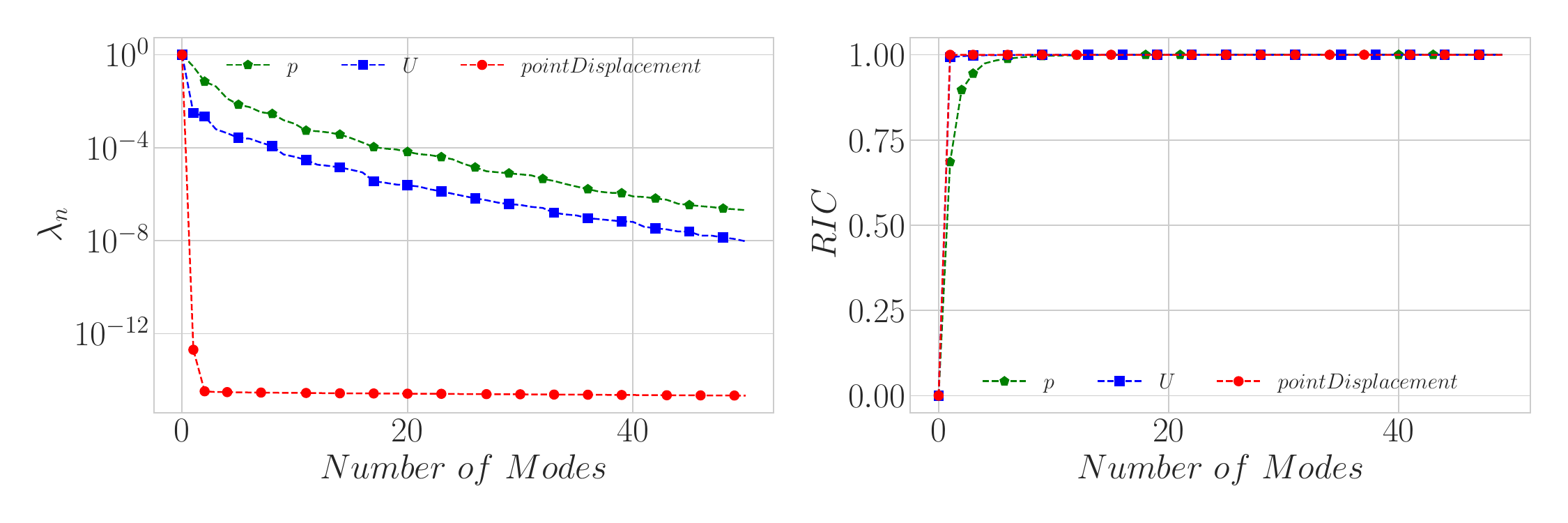} 
\caption{ Decay of eigenvalues (left) and cumulative sum of eigenvalues (right) for the POD modes. Blue lines correspond to velocity modes, green lines to pressure modes, and red lines to structural displacement modes. These plots illustrate the relative energy content captured by each mode and support the selection of truncated modal bases for reduced-order modeling.}
\label{fig:eigdecay}
\end{figure}

Fig. \ref{fig:eigdecay} illustrates the eigenvalue decay and the Relative Information Content (RIC) for the correlation matrices of the fluid velocity $\boldsymbol{u}$, pressure $p$, and  $\boldsymbol{d}^g$ fields.
The RIC, defined in \cref{RIC}, quantifies the {\RB {Kolmogorov n-width}} of the system \cite{Ahmed2020}, providing a measure of reducibility--i.e., how well a linear superposition of POD modes can approximate the full-order dynamics.
\begin{align}\label{RIC}
    RIC(M) = \left(\displaystyle\sum_{i=1}^{M}\lambda_i / \displaystyle\sum_{i=1}^{N_{s}}\lambda_i \right)\times 100,
\end{align}
where $M$ is the number of POD modes used, and $N_s$ is the total number of modes computed. RIC can then be seen as the amount of total system's energy retained by the first $M$ POD modes.

Key observations can be drawn from Fig. \ref{fig:eigdecay}.
\begin{itemize}
    \item Grid displacement $\boldsymbol{d}^g$
    \begin{itemize}
        \item Exhibits an extremely rapid eigenvalue decay, with most energy concentrated in the first POD mode (Fig. \ref{fig:eigdecay}, red line).
        \item This suggests that the Arbitrary Lagrangian-Eulerian (ALE) field can be accurately reduced to just one mode, collapsing the original 24,440 degrees of freedom (DOF) to a single DOF. Similar findings were reported in \cite{gupta2022hybrid}.
        \item The favorable reducibility likely stems from the Slerp interpolation \cite{ofsite}, which localizes displacements near the moving interface, propagating linearly toward far-field boundaries—a scenario well-suited for POD’s linear approximation.
    \end{itemize}
    \item Pressure and velocity $p$, $\boldsymbol{u}$
    \begin{itemize}
        \item Show significantly slower eigenvalue decay (Fig. \ref{fig:eigdecay}, blue/green lines), necessitating more modes for accurate ROM representation.
        \item This behavior, noted in \cite{Anttonen2003}, may arise from grid deformation effects. To mitigate this, prior work suggests:
        \begin{itemize}
            \item Domain-filtered POD to steepen eigenvalue decay, or
            \item Hadamard formulation \cite{bourguet2011reduced}, which maps deformations to a reference domain.
        \end{itemize}
    \end{itemize}
\end{itemize}

\subsubsection{ROM solution error}
Once the reconstruction error has been characterized, we aim to analyze the quality of the online problem solution. Thus, to evaluate how close the predicted ROM solutions are with respect to the FOM ones,  Fig. \ref{fig:qoi} illustrates a qualitative comparison between the solution fields contour plots corresponding to time $t = 20\ $s obtained with both the FOM and ROM solvers. The plots confirm that, to the eyeball test, the ROM solutions obtained using the mixed POD-Galerkin projection (for the fluid dynamic variables) and POD-RBF (for the grid displacement field) appear similar to the high-fidelity ones.
It is also worth pointing out that the top plots in Fig. \ref{fig:qoi} confirm the ROM is able to reproduce the 2S mode of the classical Von K\'{a}rm\'{a}n vortex street as observed in the FOM solution. 
\begingroup
\setlength{\tabcolsep}{0.05pt} 
\renewcommand{\arraystretch}{0.5} 
\begin{figure}[h!]
\centering
\begin{tabular}{ccc}
    \includegraphics[width=0.33\textwidth]{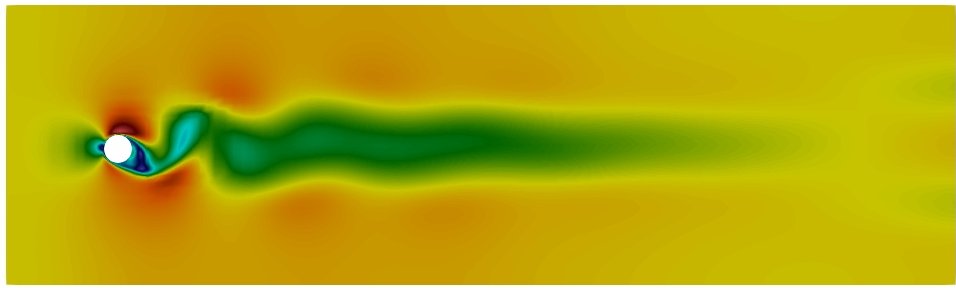} & \includegraphics[width=0.33\textwidth]{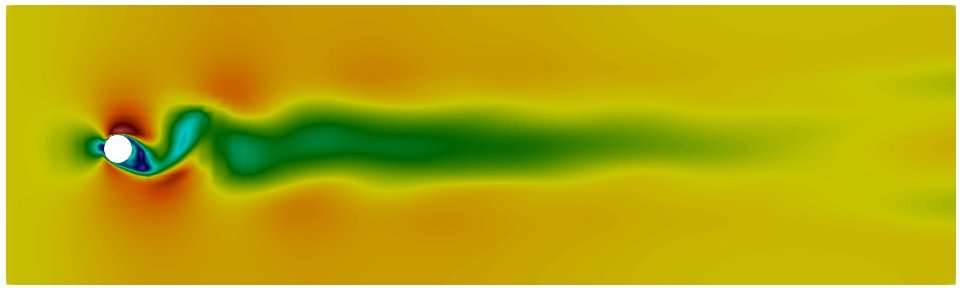}  & \includegraphics[width=0.33\textwidth]{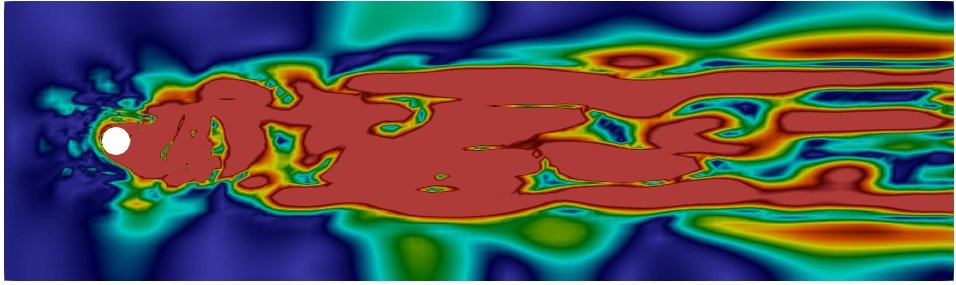}\\
    \includegraphics[width=0.23\textwidth]{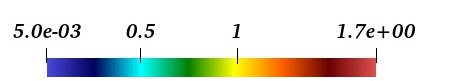} & \includegraphics[width=0.23\textwidth]{images2/U_scale.jpeg}  & \includegraphics[width=0.23\textwidth]{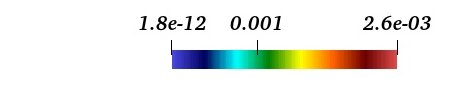}\\
    \includegraphics[width=0.33\textwidth]{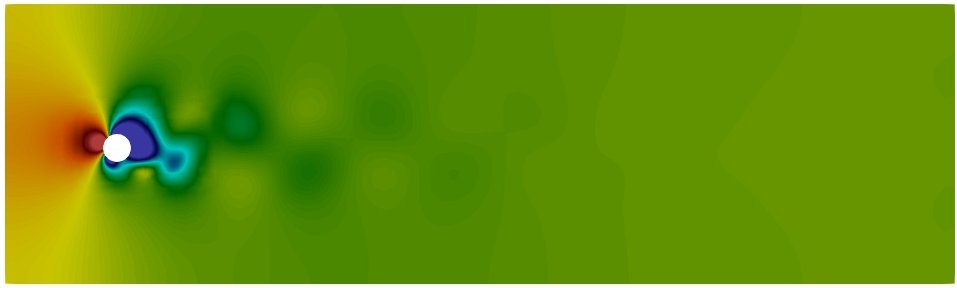} & \includegraphics[width=0.33\textwidth]{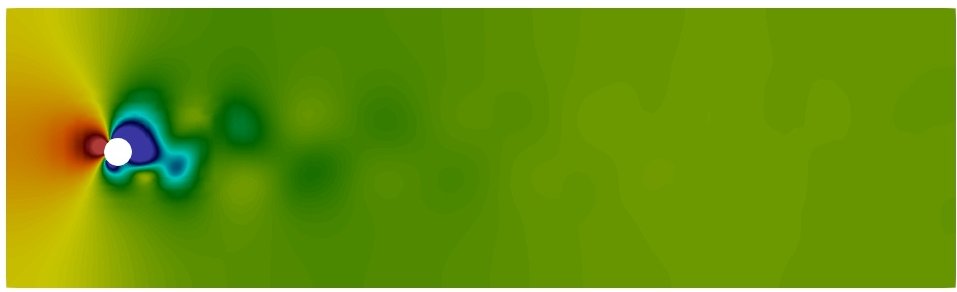}  & \includegraphics[width=0.33\textwidth]{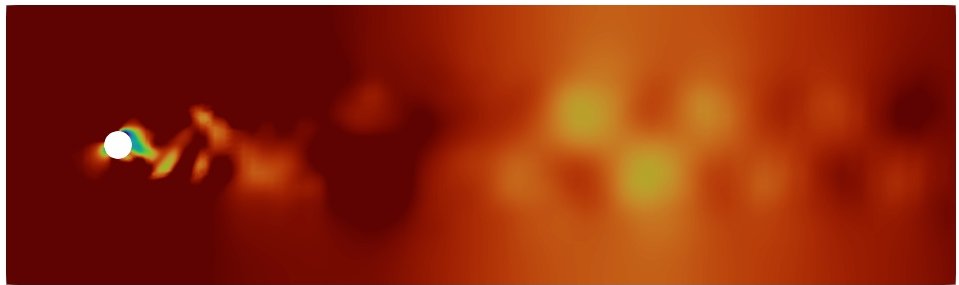}\\
    \includegraphics[width=0.23\textwidth]{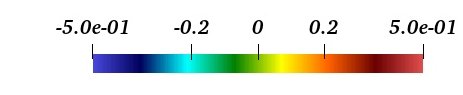} & \includegraphics[width=0.23\textwidth]{images2/p_scale.jpeg}  & \includegraphics[width=0.23\textwidth]{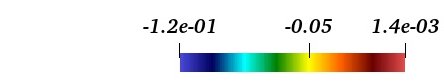}\\
    \includegraphics[width=0.33\textwidth]{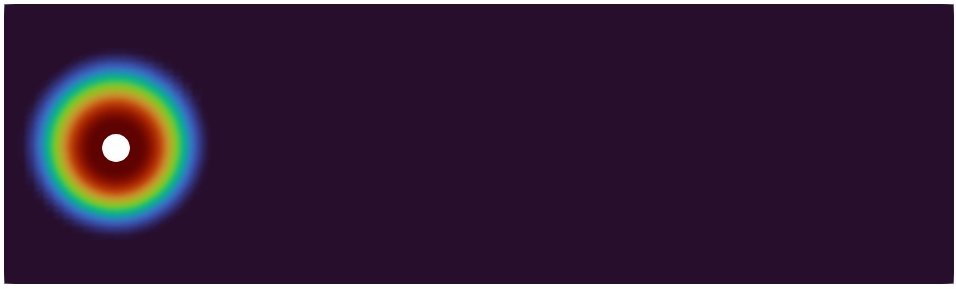} & \includegraphics[width=0.33\textwidth]{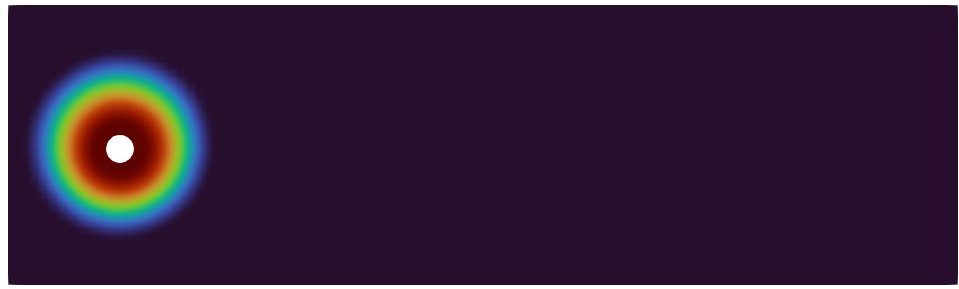}  & \includegraphics[width=0.33\textwidth]{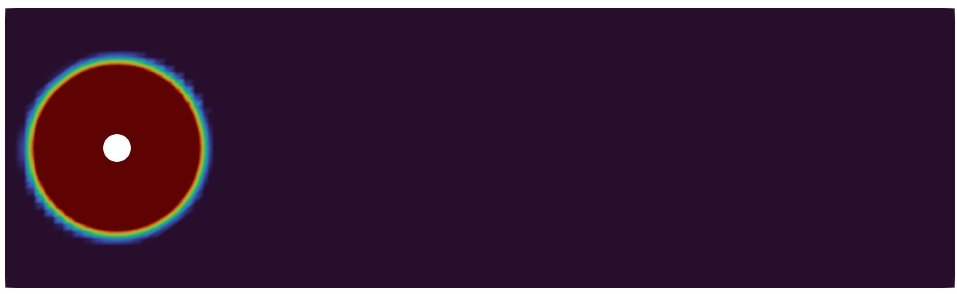}\\
    \includegraphics[width=0.23\textwidth]{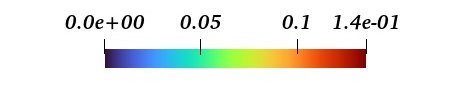} & \includegraphics[width=0.23\textwidth]{images2/displacement_scale.jpeg}  & \includegraphics[width=0.23\textwidth]{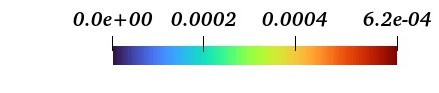}\\
\end{tabular}    
    \caption{Comparison at $t = 20\ $s between the full-order solution (first column), the reduced-order solution (second column), and the corresponding relative error maps (third column). The first row shows the velocity field, the second row the pressure field, and the third row the grid displacement. 
    The reduced-order model employs 20 POD modes for both velocity and pressure, and 3 modes for the grid displacement.}
\label{fig:qoi}
\end{figure}
\endgroup

\begin{figure}[h!]
\centering
    \includegraphics[width=\textwidth]{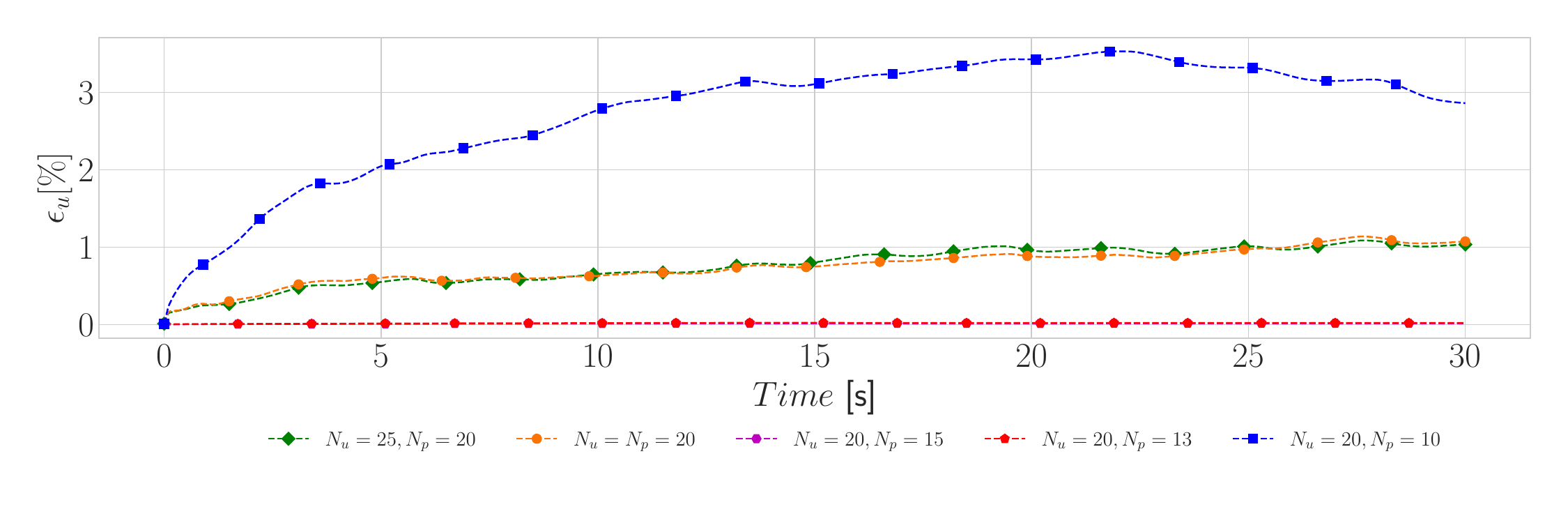}
    \caption{Velocity field relative ROM error in the $L^2$ norm as a function of time.}
    \label{fig:VelocityErrors}
\end{figure}

\begin{figure}[h!]
\centering
     \includegraphics[width=\textwidth]{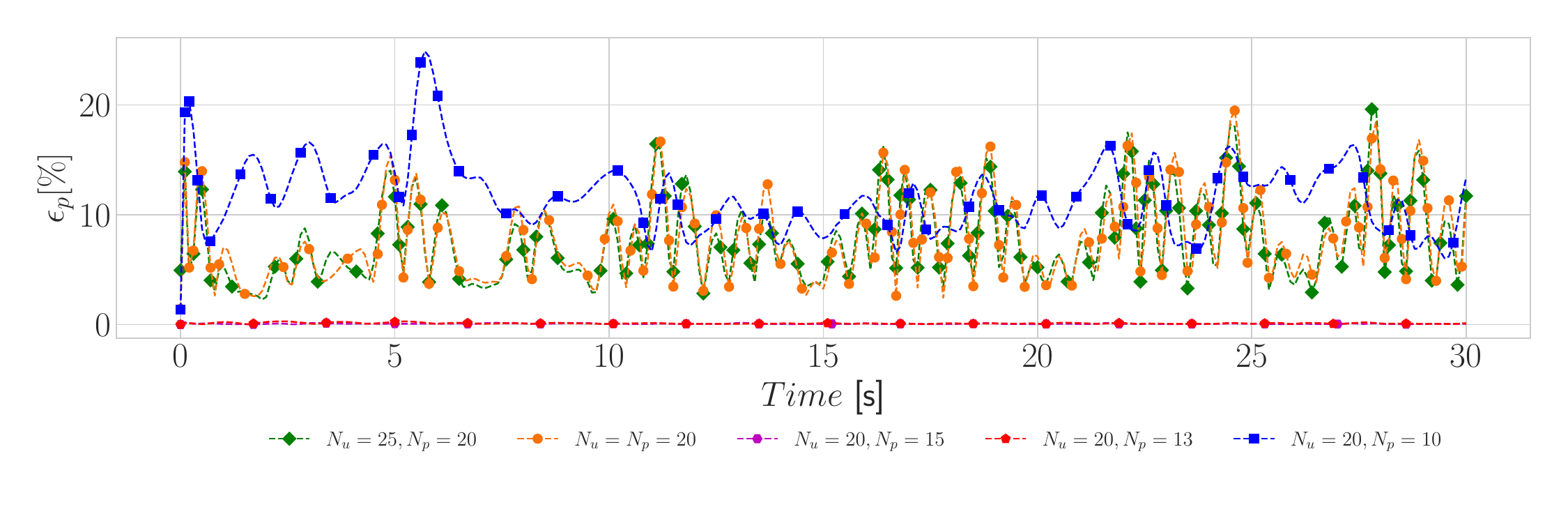}
    \caption{Pressure field relative ROM error in the L2 norm as a function of time.}
    \label{fig:PressureErrors}
\end{figure}

A more quantitative assessment of the ROM accuracy is presented in Fig. \ref{fig:VelocityErrors}, and Fig. \ref{fig:PressureErrors} which depict the time evolution of the ROM's $L^2$ relative error of the velocity and pressure fields, respectively. In the diagrams, each curve is obtained with different combination of modal truncation orders.  
The  $L^2$ relative error reported in the plots is computed, for a given quantity $q$, as
\begin{equation}
\epsilon_q = \frac{\|q_{FOM} - q_{ROM} \|_{L^2(\Omega)}}{\|q_{FOM}\|_{L^2(\Omega)}}.
\end{equation}
Figures \ref{fig:VelocityErrors} and \ref{fig:PressureErrors} illustrate the temporal evolution of velocity and pressure reconstruction errors, respectively. Both error types exhibit initial fluctuations and gradually stabilize as time progresses toward 30 seconds. This behavior suggests that the underlying numerical method is transient in nature and that the system is likely approaching a steady state.

From Fig. \ref{fig:VelocityErrors}, it is evident that the velocity error remains relatively low across all tested configurations, with a maximum around 3\%. Increasing the number of velocity POD modes ($N_u$, e.g., $N_u=25$) slightly improves accuracy, but even configurations with fewer modes (e.g., $N_u=20$) yield acceptable results. Moreover, variations in the number of pressure modes $N_p$ have only a minor effect on the velocity error, though a general trend can still be observed: higher $N_p$ values contribute marginally to improved velocity accuracy.

The configuration $N_u = N_p = 20$ provides a balanced trade-off between accuracy and computational cost, delivering low reconstruction errors in both fields. This performance aligns with expectations for POD-based reduced-order models, where velocity fields—being smoother and containing more energy—can often be effectively represented using fewer modes. This explains the relatively low velocity errors even when using only 10–15 velocity modes.

\begin{figure}[h!]
\centering
     \includegraphics[width=\textwidth]{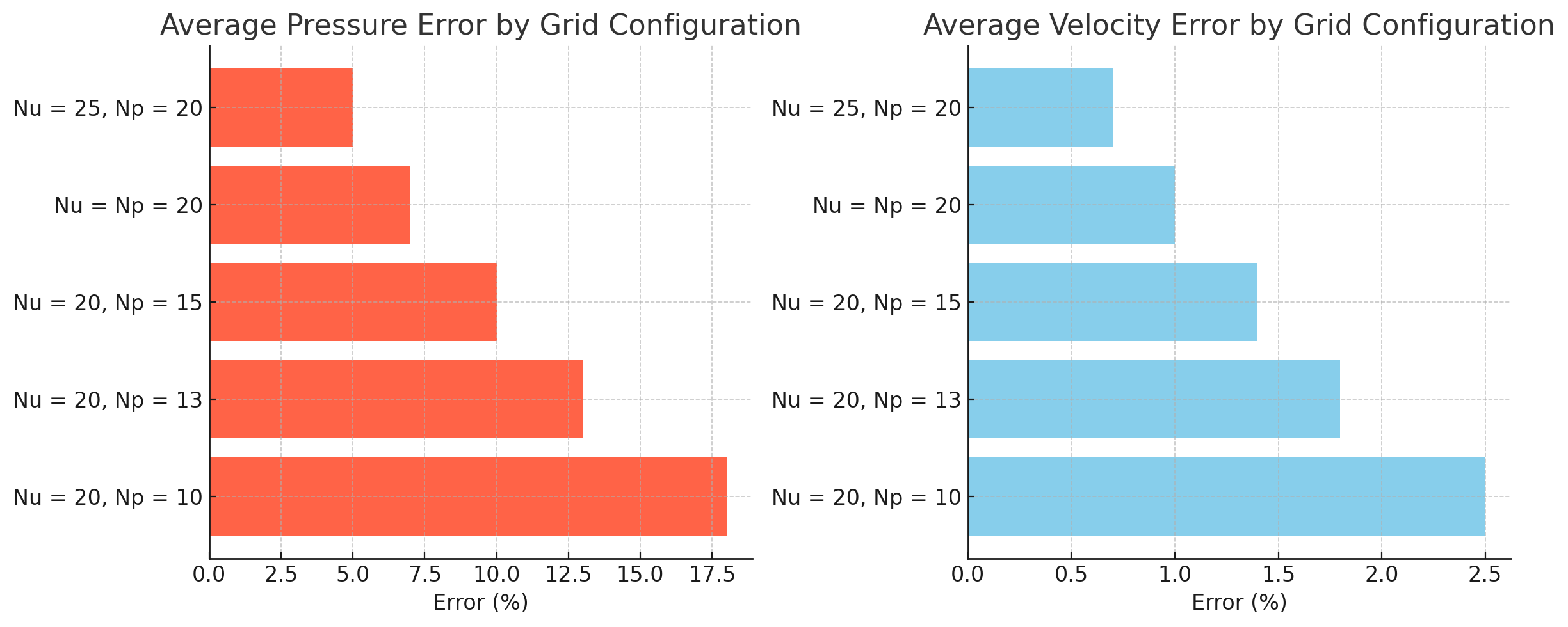}
    \caption{Summary's results of Fig. \ref{fig:PressureErrors} and Fig. \ref{fig:VelocityErrors} graphical form.}
    \label{fig:summary}
\end{figure}

In contrast, Fig. \ref{fig:PressureErrors} shows that pressure reconstruction is significantly more sensitive to the number of modes. The largest errors occur at the lowest $N_p$ (e.g., $N_p=10$), and the error decreases markedly as $N_p$ increases. This steep decline indicates that pressure requires more modes to be accurately reconstructed. This is consistent with the nature of pressure fields in incompressible flows, which often exhibit complex, global features. Moreover, pressure is typically reconstructed indirectly in projection-based ROMs (e.g., via Pressure Poisson Equation or divergence-free constraints), making it more vulnerable to errors from modal truncation.
Even at a fixed $N_u=20$, reducing $N_p$ leads to notably higher pressure errors, confirming that pressure accuracy is particularly sensitive to the number of pressure modes. 
Fig. \ref{fig:summary} summarizes the results of Fig. \ref{fig:PressureErrors} and Fig. \ref{fig:VelocityErrors} in graphical form.
While the global or time-averaged errors in the pressure and velocity fields, as shown in Figures \ref{fig:VelocityErrors} and \ref{fig:PressureErrors},  provide useful insight into the overall accuracy of the ROM, they may not fully capture \textit{localized inacuracies} that can be critical in many applications. Specifically, even a low global error may mask {high local errors} in regions of strong gradients or complex flow behavior---such as those occuring in the immediate vicinity of a body, like the cylinder in the present cross-flow simulation.

For researchers and engineers dealing with fluid dynamics problems, one of the primary objectives is often the {accuracies evaluation of the fluid dynamic forces} acting on immersed bodies or boundary surfaces. These forces are strongly dependent on the {local distribution} of the pressure and velocity near the surface. Hence, local inaccuracies---particulary around the body of interest---can significantly affect the predicted {pressure draft, lift, and shear forces}. Consequently, the global error metrics presented so far may not be sufficient indicators of how well the ROM captures the relevant physics of fluid-structure interaction studies. Inaccurate pressure and velocity predictions near cylinder surface, even if minor in a global sense, can lead to substantial errors in the {computed surface forces}, and consequently in the predicted cylinder displacement.
Therefore, a crucial next step in the ROM evaluation should include a {targeted analysis of fluid dynamic force accuracy and its influence on structural response}. This envolves computing and comparing quantities such as
\begin{itemize}
    \item Time history of lift and drag coefficients
    \item Cylinder displacement and vibration characteristics
    \item Integrated pressure and viscous forces over the body surface
\end{itemize}
Such local and integral performance metrics are indispensable for assessing the {practical applicability of ROMs in FSI simulations}, where the accuracy of force predicitions directly impacts structural safety and design decisions.
\begin{center}
\begin{figure}[h!]
    \includegraphics[width=\textwidth]{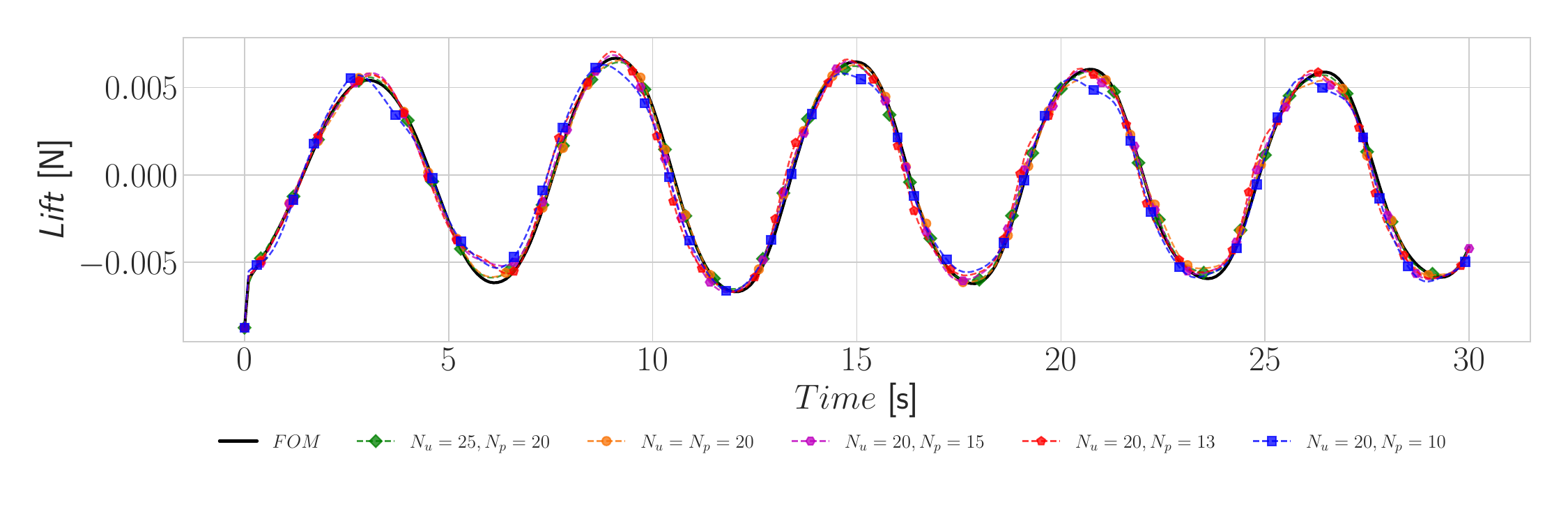}
    \caption{Time series comparison between the reference curve of the lift force acting on the cylinder in Newton unit with predicted curves  }
    \label{fig:lift}
\end{figure}
\end{center}
Fig. \ref{fig:lift} presents the time history of the lift force exerted by the fluid on the cylinder. The FOM solution is shown alongside several ROM results obtained using different combinations of velocity and pressure modes. The plot clearly demonstrates that the ROM predictions converge toward the FOM solution as the number of modes used in the online stage increases.
Notably, when 20 modes are employed for both velocity and pressure, the ROM provides a qualitatively accurate approximation of the lift force throughout the entire simulation time window. 
This observation underscores the ability of the proposed ROM methodology to reproduce key force-related dynamics with a relatively compact modal basis.
\begin{figure}[h!]
\centering
     \includegraphics[width=\textwidth]{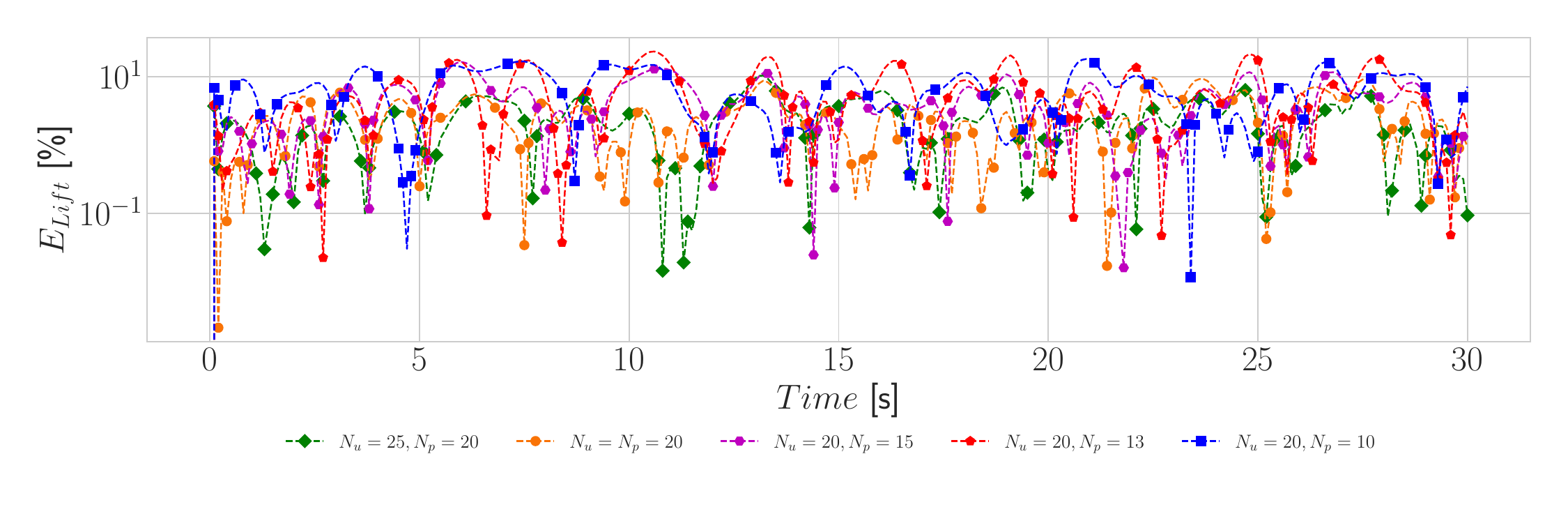}
    \caption{Time evolution of the relative errors of the pressure reduced approximation. The error values in both graphs are in percentages.}
    \label{fig:LiftForcesErrors}
\end{figure}
Further confirmation is provided in Fig. \ref{fig:LiftForcesErrors}, which shows the relative error in the lift force for the various mode combinations. As observed, when more than 20 modes are used for each field, the relative error consistently remains around or below 1\%, demonstrating the high fidelity of the ROM in capturing the fluid dynamic forces that are critical in fluid–structure interaction problems.
\begin{figure}[h!]
\centering
\includegraphics[width=\textwidth]{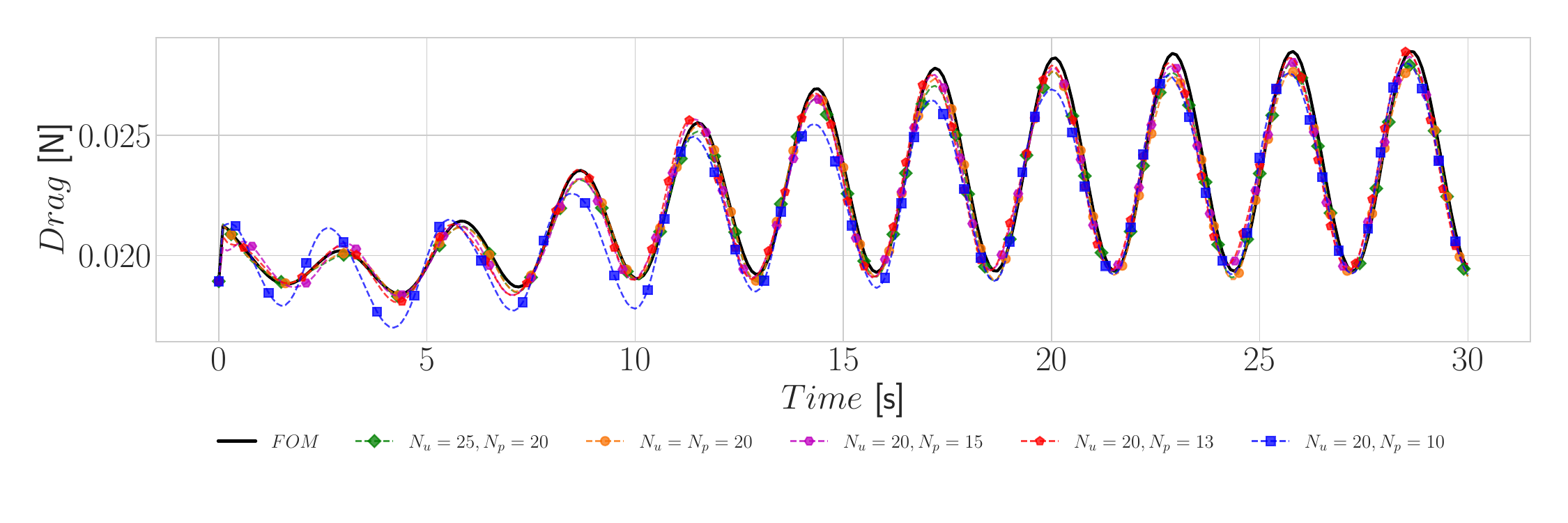} 
\caption{Time series comparison between the reference curve of the drag force acting on the cylinder in Newton unit with predicted curves. }
\label{fig:drag}
\end{figure}
Similar plots relative to the drag are presented in Fig. \ref{fig:drag} and Fig. \ref{fig:errordrag}. 
In this case, the diagram presents a comparison between the FOM drag curve and the corresponding curves obtained with ROMs that make use of different modal truncation orders. 
These plots suggest that the qualitative behavior of the cylinder resistance is well captured across all time steps of the flow simulation with a higher number of modes, and that higher absolute error appears when the number of modes decreases. 
Thus, it can be said that the accuracy shown by these plots is quite satisfactory when a number of both pressure and velocity modes higher than 20 is used. 
Additional confirmation to complement the accuracy of the ROM is shown by comparing the power spectral density curves as depicted in Fig. \ref{fig:LiftDragforcesAndPsd}.
\begin{figure}[h!]
\centering
\includegraphics[width=\textwidth]{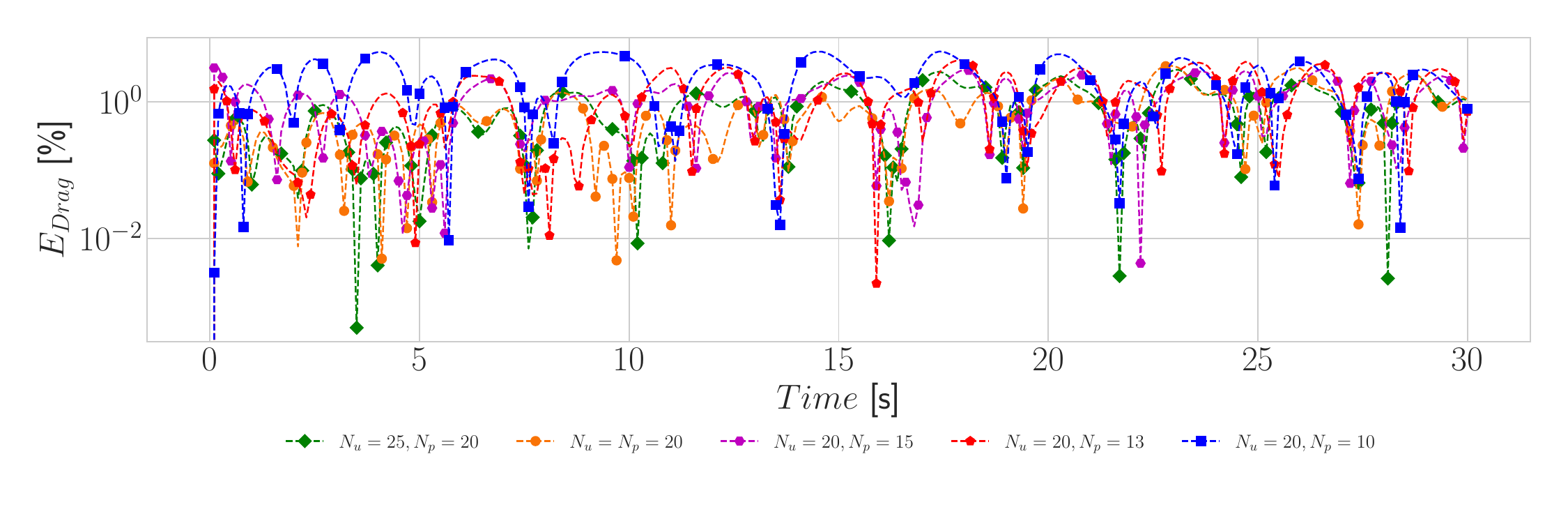}
\caption{Time series of the relative error analysis of the drag force (original and predicted signals) from Fig. \ref{fig:drag}.}
\label{fig:errordrag}
\end{figure}

The time history of the cylinder displacement provides a particularly valuable metric for assessing ROM accuracy, especially because the motion is not known a priori, unlike in cases of prescribed or forced vibrations. As such, the cylinder motion emerges naturally from the fluid–structure interaction, making it a sensitive indicator of the ROM’s predictive capabilities.

\begin{figure}[h!]
\centering
    \includegraphics[width=\textwidth]{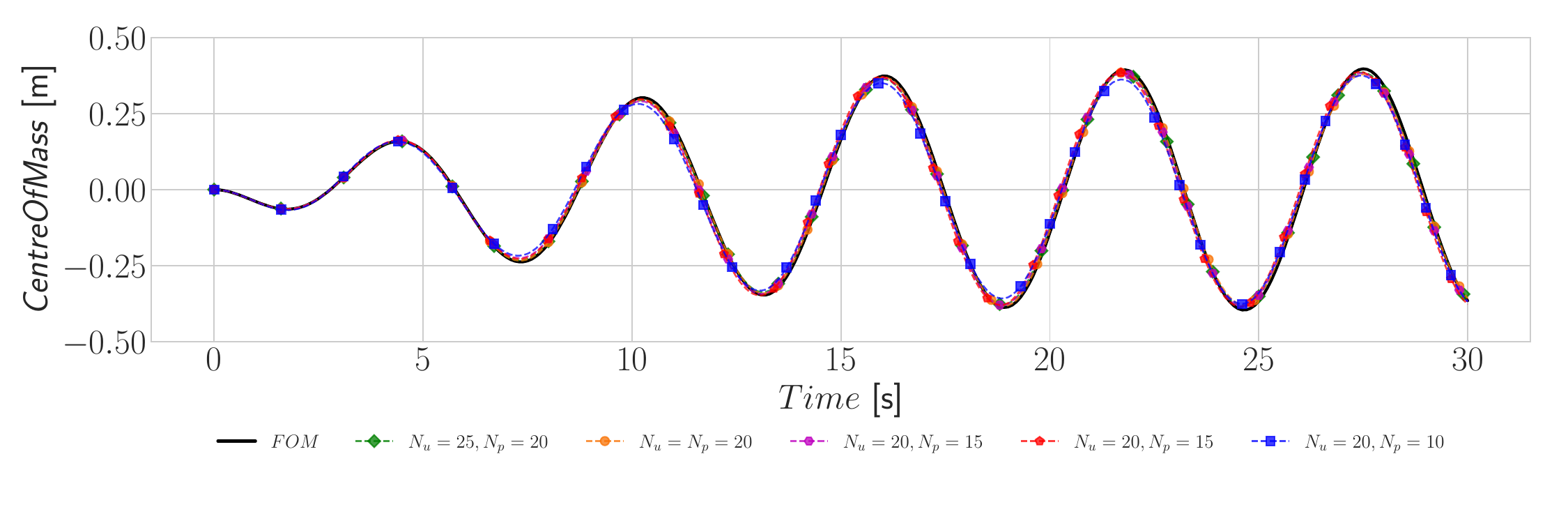}
    \caption{Time series evolution of the centre of mass.}
     \label{fig:yDisplacement}
\end{figure}

As a final verification of the proposed ROM’s accuracy, it is essential to analyze the displacement of the cylinder's center of mass over time, as computed during the simulations. Fig. \ref{fig:yDisplacement} presents a comparison between the reference high-fidelity solution and several ROM-predicted displacement curves corresponding to different modal truncation orders.

\begin{figure}[h!]
\centering
    \includegraphics[width=\textwidth]{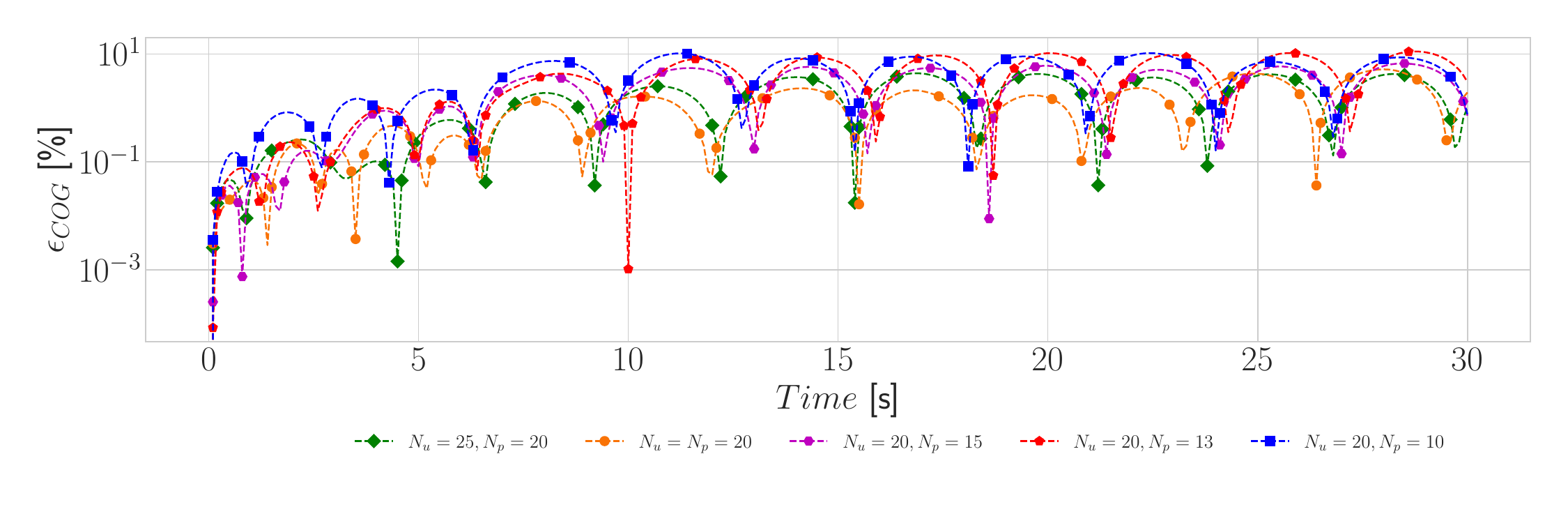}
    \caption{Time series evolution of the centre of mass relative error.}
     \label{fig:cogErrors}
\end{figure}

The results clearly show that the ROM's ability to replicate the cylinder's motion depends strongly on the number of POD modes used for velocity and pressure. 
This is expected, since the cylinder displacement is governed by the lift force, as shown in \cref{ode}, and the lift force itself is computed from the local pressure and velocity fields around the body.
In particular, the reduced solution obtained using $N_u=20$ velocity modes and $N_p=10$ pressure modes (blue line) visibly deviates from the reference curve, indicating reduced accuracy. In contrast, ROMs with higher modal resolutions exhibit much closer agreement with the reference solution. 
This is further supported by the error plots in Fig. \ref{fig:cogErrors}, which quantify the cylinder center-of-gravity error over time. For ROMs with sufficiently high modal counts, the error remains consistently below 2\% across the entire time series.
The logarithmic scale reveals differences across several orders of magnitude.
For most configurations, the relative error remains below 2\% throughout the 30-second window,
especially for $(N_u,N_p)=(25,20), (25,20), (20,20), \text{and}~ (20,15)$. The error curve for $N_p=10$ stands out with significantly higher errors, especially during the initial 10–15 seconds, reaching over 10\%.

In the implementation of POD-based ROM, the infinite-dimensional evolution equations—such as the Navier–Stokes equations—are projected onto a finite-dimensional empirical subspace, often of significantly lower dimension. A natural and fundamental question that arises in this context is: To what extent do the truncation and projection procedures preserve the long-term dynamical behavior of the original system, particularly its attractors, as discussed in \cite{Holmes1997}?
Although not shown here, additional plots confirm that the ROM is capable of accurately capturing the limit cycles observed in the FOM. This indicates that the ROM not only replicates the transient dynamics with reasonable fidelity but also faithfully approximates the attractor structure of the original high-dimensional system, thereby preserving its essential long-term behavior.

\subsubsection{Parametric sensitivity study}
This section evaluates the performance of the ROM by varying the damping coefficient $c$, a critical parameter that governs the structural response in vortex-induced vibration (VIV) systems. The damping coefficient directly influences energy dissipation and, consequently, affects both the amplitude and frequency of oscillations. Accurate modeling of these effects is essential to ensure the predictive accuracy of the ROM.

\begin{figure}[h!]
\centering
    \includegraphics[width=\textwidth]{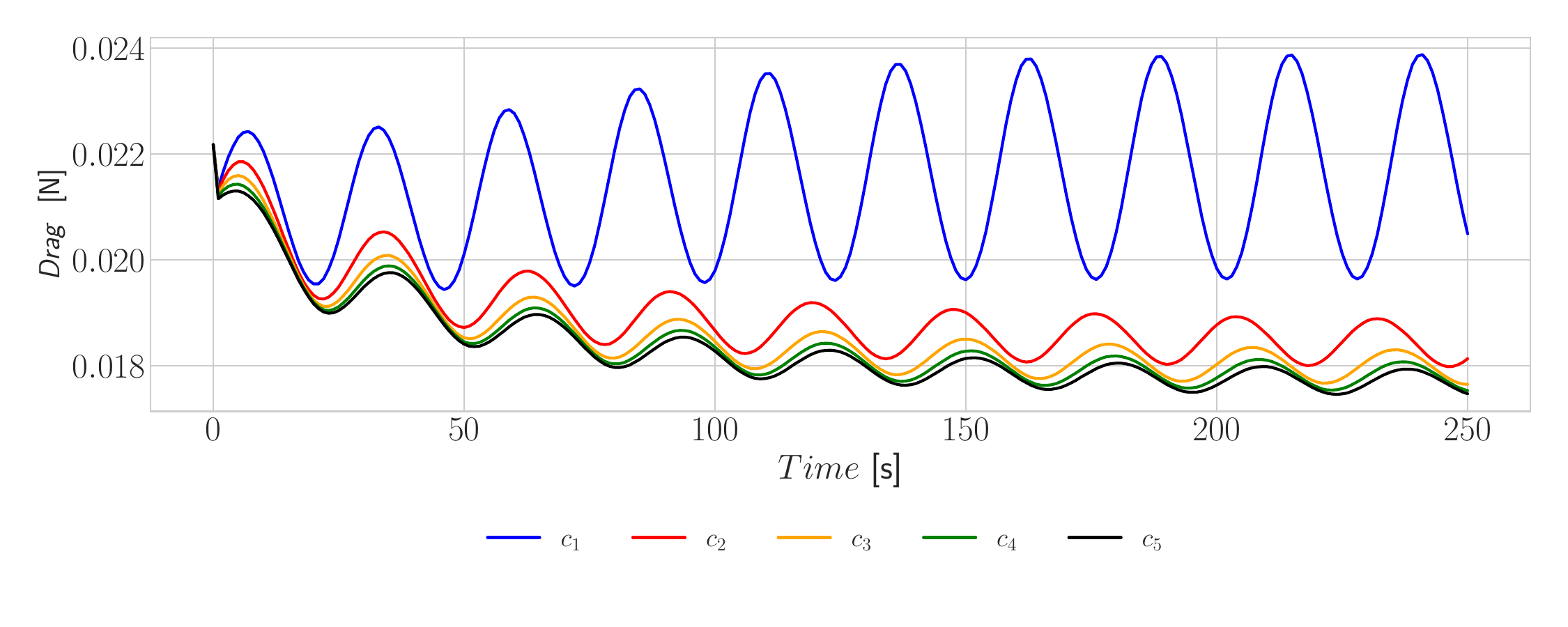}
    \caption{Time series evolution of the drag forces generated from $\mathcal{C}_{\text{training}}$ set.}
     \label{fig:training_drag}
\end{figure}

\begin{figure}[h!]
\centering
    \includegraphics[width=\textwidth]{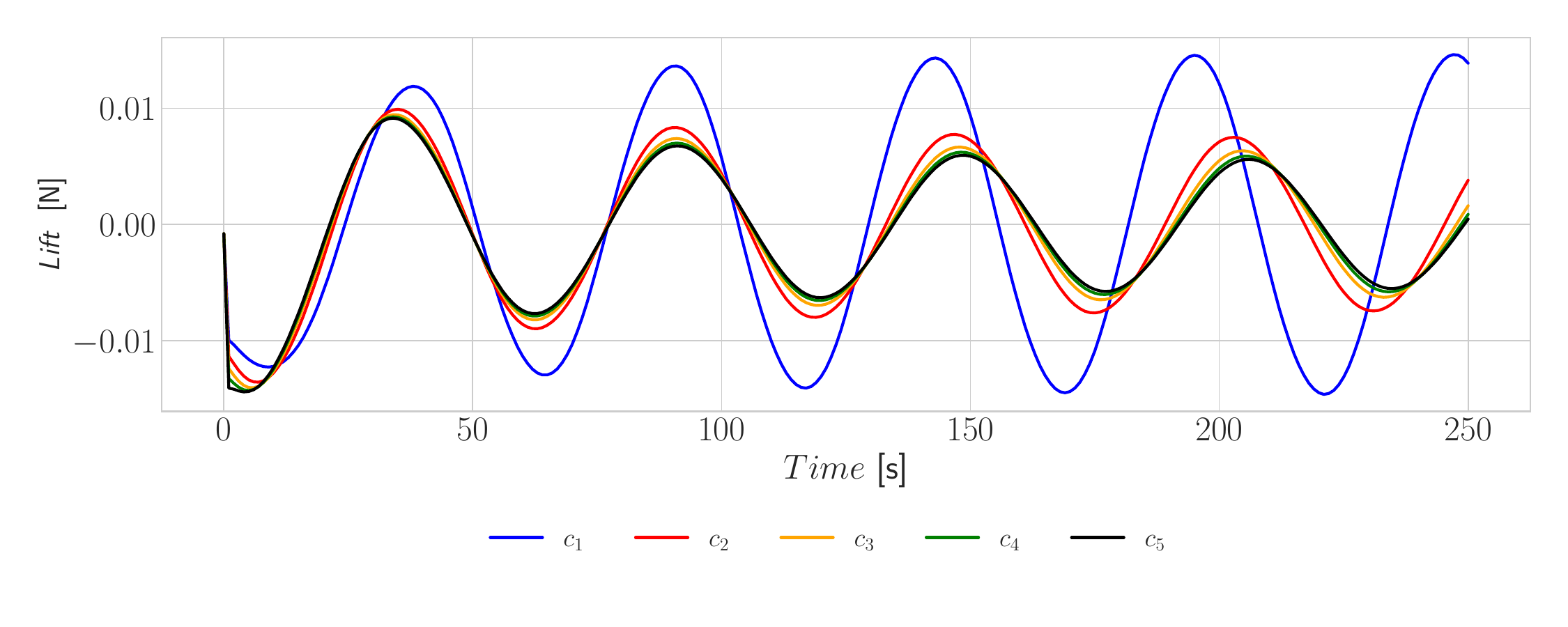}
    \caption{Time series evolution of the lift forces generated from $\mathcal{C}_{\text{training}}$ set.}
     \label{fig:training_Lift}
\end{figure}

\begin{figure}[h!]
\centering
    \includegraphics[width=\textwidth]{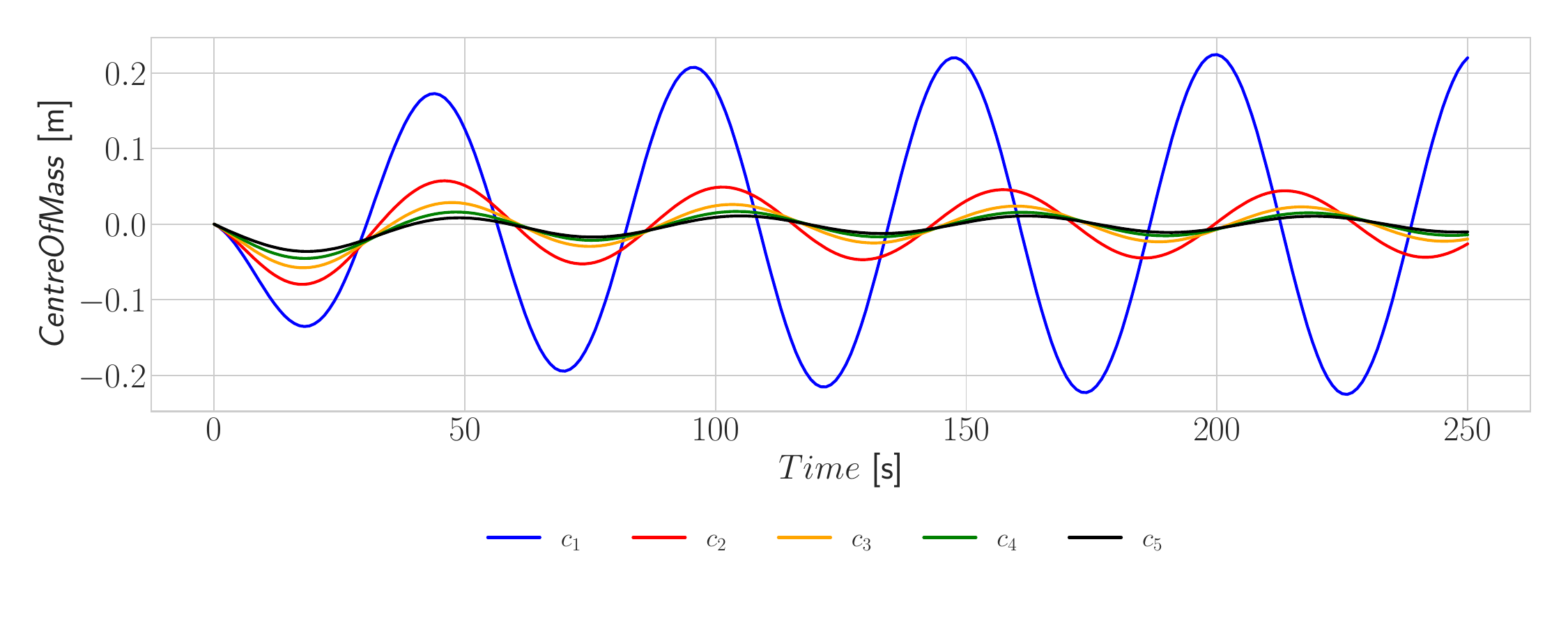}
    \caption{Time series evolution of the displacement of the centre of the mass generated from $\mathcal{C}_{training}$ set.}
     \label{fig:training_yDisplacement}
\end{figure}

To develop a ROM solver capable of capturing the dynamics of the VIV problem with a dynamic mesh, careful design of the training and testing datasets is required. The critical damping coefficient is calculated as  $c_{crit} = 2\sqrt{mk}$, $c_{crit} = 0.1644$, using  $m = 0.05kg$, and $k = 0.135115N/m$. 
The training range is defined as $c \in [c_{min}, c_{max}]$, where $c_{min} = 0.01c_{crit}$ and $c_{max} = 3c_{crit}$. 
Latin Hypercube Sampling (LHS) is employed to generate a well-distributed training set, yielding $\mathcal{C}_{training} = \{ 0.02975911, 0.13973352, 0.24222458, 0.34619364, 0.47300364 \}$. This corresponds to a damping ratio dataset $\zeta_{training} = \{ 0.1811,0.8503,1.4732,2.1062,2.8774 \}$.
Fig. \ref{fig:training_drag}, Fig. \ref{fig:training_Lift}, and 
Fig. \ref{fig:training_yDisplacement}  show the plot of the behaviors of the drag and lift forces, and  displacement of the centre of the mass when the damping coefficient vary. This set spans undamped ($c < c_{crit}$), near-critical ($c\simeq  c_{crit}$), and overdamped ($c > c_{crit}$) regimes.
It is important to ensure that the testing set values lie sufficiently close in parameter space to those used in the offline training stage to guarantee accurate ROM predictions \cite{hijazi2020data}. Based on this criterion, the testing set is chosen as $\mathcal{C}_{test} = \{0.03023731, 0.24422458, 0.50300364\}$ with the corresponding $\zeta_{test} = \{0.1839, 1.4858, 3.0604\}$.

\textbf{Training stage: }
As observed in the drag plot (Fig. \ref{fig:training_drag}), the drag force exhibits periodic oscillations across all cases. The amplitude of these oscillations decreases progressively from case $c_1$ (blue), corresponding to the first parameter in the $\mathcal{C}_{\text{training}}$ dataset, to case $c_5$ (black), which corresponds to the largest parameter in the dataset. Furthermore, the average drag force also decreases from $c_1$ to $c_5$. Case $c_1$ exhibits the highest drag values and the most pronounced oscillations, while case $c_5$ demonstrates the lowest drag and the least oscillatory behavior.
This trend suggests that as we move from case $c_1$ to $c_5$, the flow becomes more stable and  efficient, possibly due to flow condition.

The time series of the lift forces shown in Fig. \ref{fig:training_Lift} exhibit periodic behavior with a consistent trend of decreasing amplitude from case $c_1$ to $c_5$. All cases are centered around a zero-mean lift, indicating symmetric oscillations and no net lift over time. Among them, case $c_1$ displays the largest peak-to-peak variation, suggesting more pronounced unsteady flow behavior.
The progressive reduction in lift fluctuations from $c_1$ to $c_5$ indicates enhanced flow stability and a decrease in unsteady fluid dynamic forces acting on the oscillating cylinder. This trend is favorable for reducing structural fatigue and improving controllability in flow-structure interaction scenarios.

Figure~\ref{fig:training_yDisplacement} shows the vertical displacement of the center of mass over time for five configurations, labeled $c_1$ to $c_5$ and color-coded from blue to black. These configurations result from varying the damping coefficients within $\mathcal{C}_{\text{training}}$. All cases exhibit periodic oscillations, indicative of a dynamic unsteady flow regime. This behavior suggests that the system experiences cyclic vertical motion, likely driven by vortex-induced vibrations (VIV).
Building on this observation, the displacement amplitudes show a clear trend: they decrease monotonically from $c_1$ (blue, largest amplitude) to $c_5$ (black, smallest amplitude). For instance, $c_1$ exhibits peak displacements exceeding $\pm0.2$ m, whereas $c_5$ remains confined within $\pm0.02$ m. This trend reflects the damping effect's growing influence in suppressing oscillations.
Despite this steady reduction in amplitude, the oscillation frequency remains nearly constant across all configurations. This indicates that the dominant frequency of the fluid forcing—presumably from vortex shedding—remains stable. Therefore, it is primarily the structural response magnitude, not the excitation frequency, that varies with damping.

Finally, it is noteworthy that the oscillations remain in phase across all configurations. This phase alignment suggests a consistent external periodic excitation, most likely due to synchronized vortex shedding. Additionally, it may point to the presence of a phase control mechanism—either active or passive—ensuring coherent motion.

\textbf{Online stage: }
To test the online performances, 3 new scalar damping coefficients have been randomly selected.  
40 modal basis functions have been picked for the prediction of velocity, 35 modal basis functions for pressure  while 3 modal basis functions have been employed for grid nodes displacement. 
The increased number of modes required for accurate representation in this framework was attributed to the
relative motion of grid points caused by deformation as it has been already discussed in \cite{anttonen2005applications,anttonen2001techniques}. 
Moreover, as deformation increased, measured by the number of moving grid points, the number of modes required for an accurate solution also increased.

\begin{figure}[h!]
\centering
    \includegraphics[width=\textwidth]{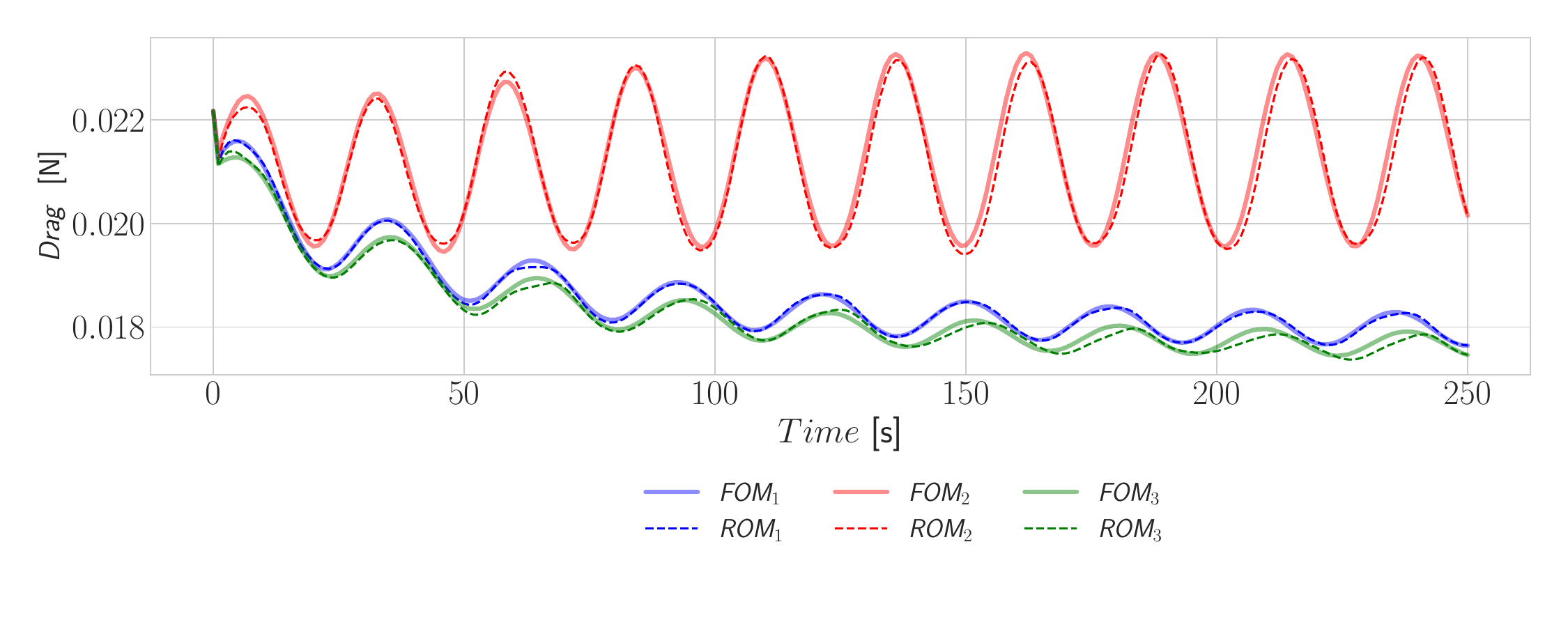}
    \caption{Comparison of drag forces time series between FOMs and ROMs predictions generated from $\mathcal{C}_{test}$ dataset.}
     \label{fig:test_drag}
\end{figure}
As shown in Fig. \ref{fig:test_drag}, the ROMs accurately capture the frequency and phase of the FOMs, although slight discrepancies in amplitude are observed. Among the models, $FOM_{1}/ROM_{1}$ exhibit the highest amplitudes, while $FOM_{3}/ROM_{3}$ show the lowest, reflecting diminished unsteady flow dynamic forces. Overall, the drag responses show strong agreement between ROM and FOM, supporting the validity of the ROMs.
\begin{figure}[h!]
\centering
    \includegraphics[width=\textwidth]{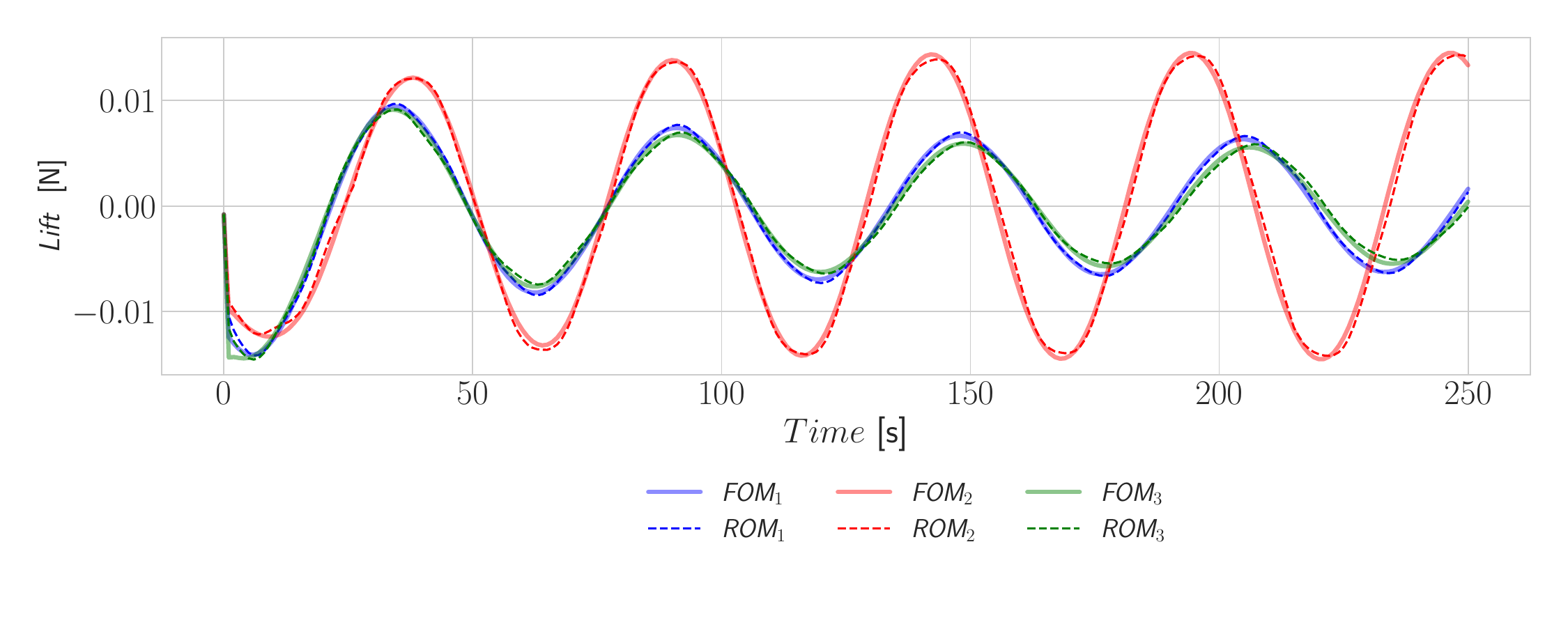}
    \caption{Comparison of lift force time series between FOMs and ROMs predictions generated using the $\mathcal{C}_{test}$ dataset}
     \label{fig:test_Lift}
\end{figure}

Similarly, visual inspection of Fig. \ref{fig:test_Lift} indicates that the ROMs closely replicate the lift behavior of the FOMs, particularly in terms of frequency and phase alignment. A consistent reduction in amplitude is observed from $FOM_{1}/ROM_{1}$ to $FOM_{3}/ROM_{3}$, suggesting increased flow stability or reduction in vortex-induced oscillations. The very low amplitude responses in $FOM_{3}$ and $ROM_{3}$ further indicate a near steady-state flow condition.
\begin{figure}[h!]
\centering
    \includegraphics[width=\textwidth]{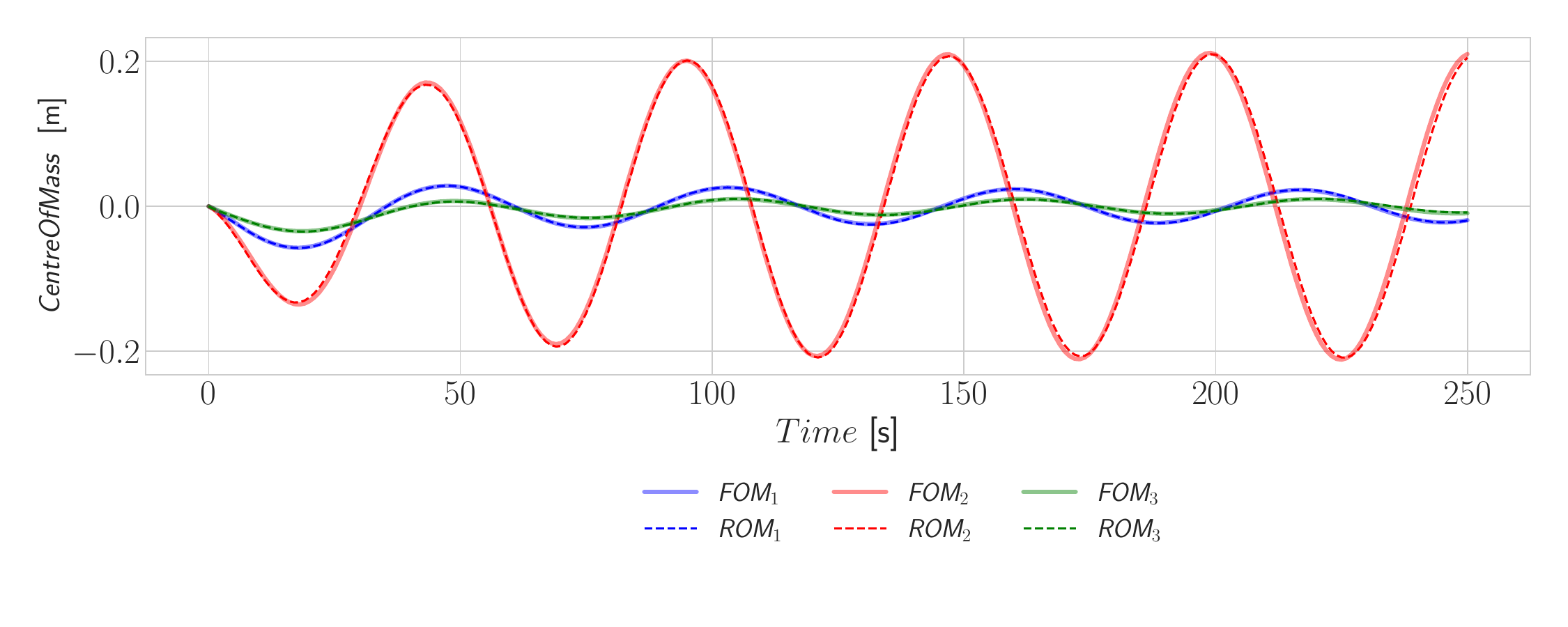}
    \caption{Comparison between the time series of the centre of the mass displacement of FOMs and ROMs prediction generated from $\mathcal{C}_{test}$ dataset.}
     \label{fig:test_yDisplacement}
\end{figure}
In Fig. \ref{fig:test_yDisplacement}, all models demonstrate periodic oscillations. The ROMs successfully capture both the frequency and overall trend of the FOM responses, with only minor amplitude deviations. Notably, $ROM_{1}$ aligns closely with $FOM_{1}$, indicating high predictive fidelity. The gradual reduction in amplitude from $FOM_{1}/ROM_{1}$ to $FOM_{3}/ROM_{3}$ suggests increasingly damped behavior or progressive flow stabilization. These findings collectively validate the ROMs’ capability to accurately predict the dynamic flow response under unseen parameter variations.
\section{Conclusion and outlooks}\label{conclusion}
This paper presents high-fidelity and surrogate simulations for a flow passing a cylinder with a moving mesh moving boundary in the Arbitrary Lagrangian-Eulerian approach. 
The computational mesh deformation is considered a part of the solution state vector while constructing the reduced POD basis. 
The method is demonstrated by using two case studies involving vortex-induced vibration of a cylinder at Reynolds number { \RB{$Re=200$} }. 
In the first study, time is treated as the sole parameter, while the second includes both time and the structural damping coefficient as parametric inputs. 
The results demonstrate that reduced-order models (ROMs) can be effectively developed for the PIMPLE algorithm with moving boundaries using the ALE formulation within the OpenFOAM framework. The constructed ROM successfully captures the essential dynamics of VIV in a laminar flow regime, in agreement with full-order simulations performed using OpenFOAM solvers.

Future work will focus on incorporating hyper-reduction strategies tailored for fluid–structure interactions involving moving boundaries and geometric parameterizations. 

\section*{Acknowledgments}

Gianluigi Rozza acknowledges funding from the European Union Horizon 2020 Program in the framework of European Research Council Executive Agency: H2020 ERC CoG 2015 AROMA-CFD project 681447 "Advanced Reduced Order Methods with Applications in Computational Fluid Dynamics" P.I. Professor Gianluigi Rozza. Gianluigi Rozza also acknowledges funding from the Italian Ministry of University and Research (MUR) in the form of PRIN ``Numerical Analysis for Full and Reduced Order Methods for Partial Differential Equations'' (NA-FROM-PDEs) project, and by INdAM GNCS.
Andrea Mola acknowledges the financial support under the National Recovery and Resilience Plan (NRRP), Mission 4, Component 2, Investment 1.1, Call for tender No. 1409, funded by the European Union – NextGenerationEU – Project Title ROMEU – CUP D53D2301888001 - Grant Assignment Decree No. 1379 adopted on 01/09/2023 by the Italian Ministry of University and Research (MUR) and the financial support by INdAM GNCS.
Giovanni Stabile acknowledges the financial support under the National Recovery and Resilience Plan (NRRP), Mission 4, Component 2, Investment 1.1, Call for tender No. 1409, funded by the European Union – NextGenerationEU– Project Title ROMEU – CUP P2022FEZS3 - Grant Assignment Decree No. 1379 adopted on 01/09/2023 by the Italian Ministry of University and Research (MUR), and acknowledges the financial support by the European Union (ERC, DANTE, GA-101115741). Views and opinions expressed are however those of the author(s) only and do not necessarily reflect those of the European Union or the European Research Council Executive Agency. Neither the European Union nor the granting authority can be held responsible for them.

\appendix
\appendix

\section{Synchronization  analysis}\label{append1}
\begin{center}
\begin{figure}[h!]
    \includegraphics[width=\textwidth]{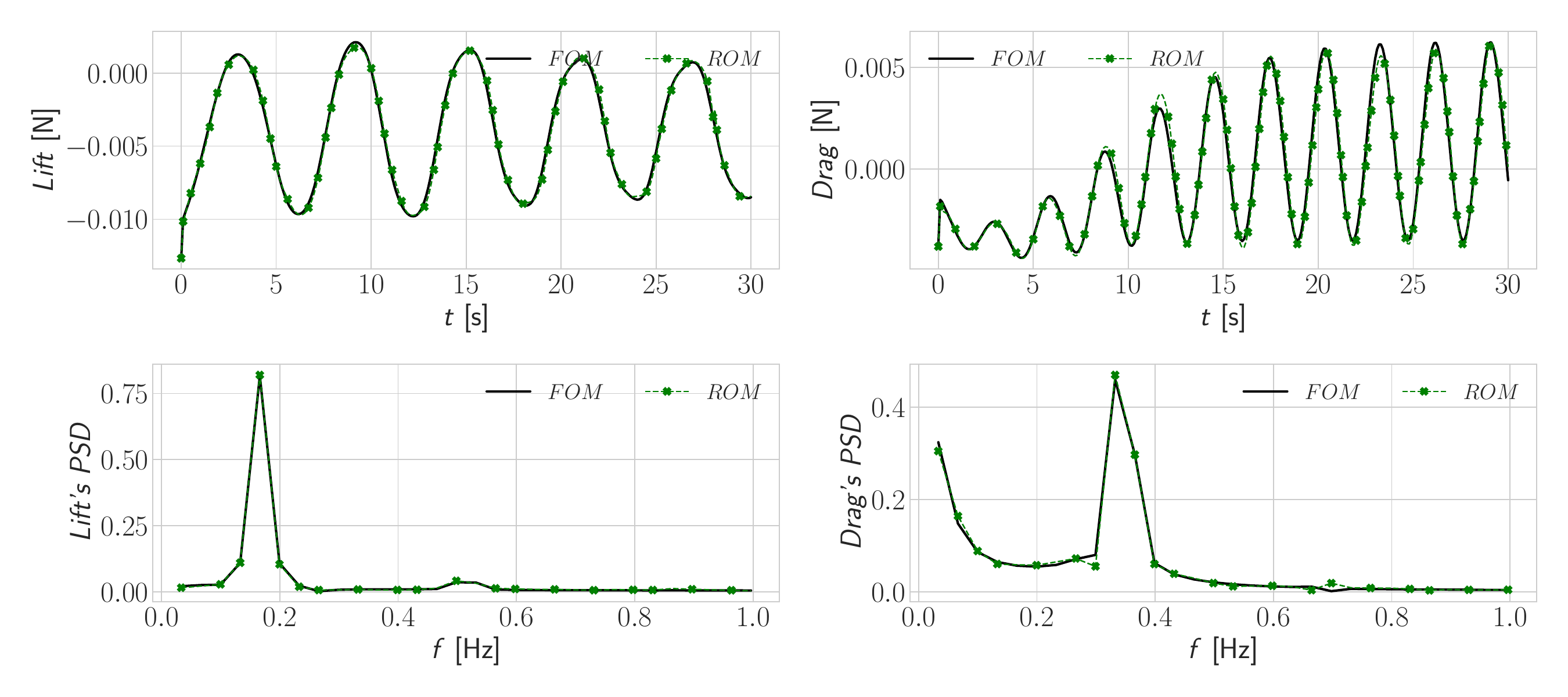}
    \caption{First row, from left to right: the time histories of the lift and drag forces.  The solid black lines are the FOM curves, and the dashed green line are the ROM curves obtained with 16 modes for the velocity and 21 modes for the pressure. Second row, from left to right: Power spectra density comparison of the lift and drag coefficients}
    \label{fig:LiftDragforcesAndPsd}
\end{figure}
\end{center}

The periodic state reached is characterized by the oscillation of the drag coefficient at twice ($f_{drag} \approx 2f_{sh}$) the lifting frequency \cite{Placzek2009} as one can see in the right plot of Fig. \ref{fig:LiftDragforcesAndPsd}.
One of the most exciting characteristics of the fluid body interaction is that of synchronization, or "lock-in," between the vortex shedding and the cylinder vibration frequencies. When the wake is synchronized, the vortex-shedding frequency diverges from that corresponding to a fixed cylinder. 
It becomes equal to the frequency of the cylinder oscillation, as shown in Fig. \ref{fig:LiftDragforcesAndPsd}.

\begin{figure}[h!]
\centering
\begin{tabular}{cc}
\includegraphics[width=0.33\textwidth]{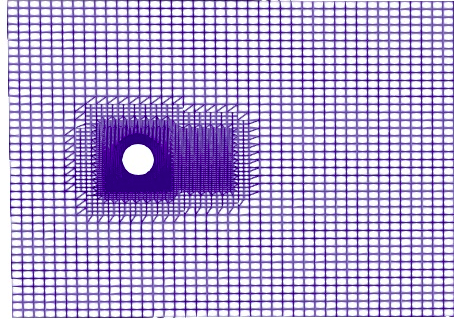}
    &
\includegraphics[width =0.35\textwidth]{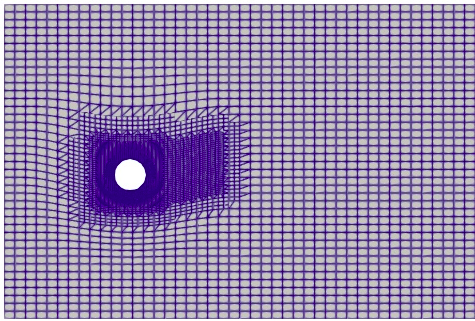}\\
     (a) & (b)
\end{tabular}    
    \caption{The  mesh used in the simulations: (a) the initial mesh, (b) the deformed mesh in correspondence with a cylinder displacement of approximately 40\% of the diameter.}
    \label{fig:mylabel}
\end{figure}

\section{Reduced-PIMPLE: Core Computational Steps}\label{reduced-pimple-steps}
\begin{enumerate}

\item \textbf{Compute the Forces}  
\begin{itemize}
    \item Based on predicted velocity \( \boldsymbol{u}_h^{n*} \) and previous pressure \( \boldsymbol{p}_h^{n-1} \).
    \item Used to evaluate the fluid-induced forces on the structure.
\end{itemize}

\item \textbf{Solve the Rigid Body Problem}  
\begin{itemize}
    \item Integrate Newton’s equations of motion to update the position of the structure (e.g., \( y_C^{\text{new}} \) for a cylinder).
\end{itemize}

\item \textbf{Evaluate RBF Map and Reconstruct Grid Motion}  
\begin{itemize}
    \item Compute RBF coefficients using updated structure position.
    \item Reconstruct displacement field: \( \boldsymbol{\delta}^n = \bm{\Psi} \boldsymbol{c}^T \).
    \item Update the computational mesh accordingly (a cheap operation in OpenFOAM).
\end{itemize}

\item \textbf{Assemble the Momentum Equation}  
\begin{itemize}
    \item Using the updated mesh, assemble the system matrix \( \mathbf{A}_u \) and RHS vector \( \mathbf{b}_u \).
\end{itemize}

\item \textbf{Solve the Reduced System and Reconstruct Velocity}  
\begin{itemize}
    \item Solve: \( \boldsymbol{\Phi}^T \mathbf{A}_u \boldsymbol{\Phi} \boldsymbol{a}^* = \boldsymbol{\Phi}^T \mathbf{b}_u \).
    \item Reconstruct velocity: \( \boldsymbol{u}_h^{n*} = \boldsymbol{\Phi} \boldsymbol{a}^* \).
\end{itemize}

\item \textbf{Assemble the Pressure Poisson Equation (PPE)}  
\begin{itemize}
    \item Construct PPE matrix using divergence of predicted velocity and continuity constraint.
\end{itemize}

\item \textbf{Solve the Reduced PPE and Reconstruct Pressure}  
\begin{itemize}
    \item Solve: \( \boldsymbol{\Xi}^T \mathbf{A}_p \boldsymbol{\Xi} \boldsymbol{b}' = \boldsymbol{\Xi}^T \boldsymbol{b}_p \).
    \item Reconstruct pressure increment: \( \boldsymbol{p}' = \boldsymbol{\Xi} \boldsymbol{b}' \).
\end{itemize}

\item \textbf{Momentum Corrector}  
\begin{itemize}
    \item Compute corrected velocity: \( \boldsymbol{u}' = \boldsymbol{u}_h^{n*} - \mathbf{A}^{-1}\mathbf{B}_p \boldsymbol{p}' \).
    \item If enabled, loop through inner pressure corrections (PISO).
\end{itemize}

\end{enumerate}

\bibliographystyle{plain}
\bibliography{biblio}
\end{document}